%AMS-TEX Ver. 2.1

\input amstex
\mag=\magstep1
\documentstyle{amsppt}

\hcorrection{3truemm}
\vcorrection{-7.2truemm}

\def\r{{\Bbb R}^d}

\def\mm#1{\left[#1 \right]}

\def\RE{\operatorname{Re}}

\def\dis{\displaystyle}

\CenteredTagsOnSplits

\topmatter
\title
Asymptotic expansion of the expected volume
of the Wiener sausage in even dimensions
\endtitle
\leftheadtext{Y. HAMANA}
\author
Yuji Hamana
\endauthor
\affil
Department of Mathematics, Kumamoto University
\endaffil
\rightheadtext{THE WIENER SAUSAGE IN EVEN DIMENSIONS}
\thanks
Partly supported by the Grant-in-Aid for Scientific Research (C)
No.24540181, Japan Society for the Promotion of Science.
\endthanks
\abstract
We consider the Wiener sausage for a Brownian motion
up to time $t$ associated with a closed ball in even dimensional cases.
We obtain the asymptotic expansion of the expected volume of
the Wiener sausage for large $t$. The result says that
the expansion has many $\log$ terms, which do not appear
in odd dimensional cases.
\endabstract
\endtopmatter

\document

%111111111111111111111111111111111111111111111111111

\head
1. Introduction
\endhead

\noindent
The Wiener sausage appears in various kinds of situations, for example
heat conduction problems, random Schr{\"o}dinger operators and
Brownian motion in random obstacles.
In particular, in connection with
heat conduction problems, the volume of the Wiener sausage
on the time interval $[0,t]$ for a Brownian motion
associated with a non-polar compact set has been investigated
for a long time. The expected volume of the Wiener sausage
is interpreted as the total energy flow from the non-polar set.
For large $t$ it is asymptotically equal to $2\pi t/\log t$
in the two dimensional case and $t$ multiple of the Newtonian capacity
of the non-polar set in higher dimensions, which can be found in \cite{17}.
The improvement of these results are in \cite{14, 15}
and the same problem for a stable sausage is discussed
in \cite{3, 16}.

We should mention that limit theorems for the volume of the Wiener sausage
have been established. The laws of large numbers are described
in \cite{10, 12}. The central limit theorems are proved in \cite{13} and
the results concerning large deviations are given in
\cite{1, 2, 7}.

This article deals with the case that the non-polar compact set is
a closed ball with radius $r$.
If the dimension is odd, the explicit form of the mean volume of
the Wiener sausage and its asymptotic expansion
have been obtained in \cite{5, 6}. They are represented by zeros
of a suitable modified Bessel function.
In comparison with the odd dimensional case,
the explicit form in the even dimensional case
has an additional part which concerns the integral of
some function consisting of modified Bessel functions.
The formula is given in \cite{9}. 
The asymptotic behavior of the mean volume of
the Wiener sausage is also discussed in \cite{9}.
It is remarkable that it has the term of $\log t/t^{d-3}$,
which does not appear in the odd dimensional case.
The purpose of this paper is to give the asymptotic expansion
of the expected volume of the Wiener sausage, which means
the improvement of  the result in \cite{9}.

This article is organized as follows. Section 2 is devoted to
discussing about the expected volume of the Wiener sausage
and giving its asymptotic behavior in the even dimensional cases.
Section 3 deals with asymptotic behavior of a function
consisting of modified Bessel functions,
which plays an important role to prove the results given
in Sections 2. We consider the case when the dimension is even and
not less than six in Section 4
and discuss the four dimensional case in Section 5. Section 6
is devoted to proving lemmas given in Section 4.

%2222222222222222222222222222222222222222222222

\head
2. The expected volume of the Wiener sausage
\endhead

\noindent
Let $A$ be a compact set in $\r$ and $\{B(t)\}_{t\geq0}$ be a Brownian motion
on $\r$. The Wiener sausage $\{W(t;A)\}_{t\geqq0}$
for the Brownian motion associated with $A$ is the process
defined by
$$
W(t;A)=\{x\in\r\,;\,x+B(s)\in A\text{ for some $s\in[0,t]$}\}
$$
for $t\geqq0$. It is easy to show that $W(t;A)$ is also compact in $\r$.
For $t>0$ let
$$
L(t;A)=\int_{\r\setminus A}P_x[\tau_A\leqq t]dx,
$$
where $\tau_A=\inf\{t\geqq0\,;\,B(t)\in A\}$
and $P_x$ is the probability measure of events related to
the Brownian motion starting from $x\in\r$.
It is easy to see that the expectation of the volume of $W(t;A)$
coincides with the sum of $L(t;A)$ and the volume of $A$ (cf. \cite{17}).

We first consider the formula for  $L(t;A)$ when $A$ is a ball.
Let $D$ be a closed ball with center $0$ and radius $r$.
The notation $L(t)$ will be used to denote $L(t;D)$.
One and three dimensional cases are easy.
Indeed, by the well-known formula for $P_x[\tau\leqq t]$,
we can see directly that
$$
L(t)=
\cases
2\sqrt{2t/\pi}\quad&\text{if $d=1$,}\\
2\pi rt+4r^2\sqrt{2\pi t}\quad&\text{if $d=3$.}
\endcases
$$
For details, see \cite{5, 10, 14}.
In higher dimensional cases, although the explicit form of
$P_x[\tau\leqq t]$ has been obtained in~\cite{8},
it seems to be difficult to carry out the integration on~$x$.
However the Laplace transform of $L$ is given by
$$
\int_0^\infty e^{-\lambda t}L(t)dt=\frac{S_{d-1}r^{d-1}}{\sqrt{2\lambda^3}}
\frac{K_{d/2}(r\sqrt{2\lambda})}{K_{d/2-1}(r\sqrt{2\lambda})}
\tag2.1
$$
for $\lambda>0$, where $S_{d-1}$ has been used for
the surface area of $d-1$ dimensional unit sphere
and $K_\nu$ is the modified Bessel function of the second kind
(cf. \cite{5}). It is shown in \cite{6, 9} that
the right hand side of (2.1) can be inverted.

Before giving the result on $L(t)$, we recall several facts
on modified Bessel functions. For each complex number $\nu$
the modified Bessel function of order $\nu$ is
the fundamental solutions of
the modified Bessel differential equation
$$
z^2\frac{d^2w}{dz^2}+z\frac{dw}{dz}-(z^2+\nu^2)w=0.
\tag2.2
$$
The standard notation $I_\nu$ and $K_\nu$ are used to denote
the solutions, which are called the first kind and
the second kind, respectively (cf. \cite{11, 18}).
When $\nu$ is real, both $K_\nu(x)$ and $I_\nu(x)$
are real and positive for $x>0$.
In this paper we need information on the zeros of $K_\nu$.
For $\nu\in\Bbb R$ let $\phi(\nu)$
be the number of zeros of $K_\nu$. It is known
that $\phi(\nu)$ is equal to $|\nu|-1/2$ if $\nu-1/2$ is an integer
and the even number closest to $|\nu|-1/2$ otherwise.
In particular, $\phi(\nu)=0$ if $|\nu|<3/2$.
$K_\nu$ has one real zero if $|\nu|=2n+3/2$ for some integer~$n$
and no zero otherwise. When $\phi(\nu)\geqq1$,
we write $z_{\nu,1},z_{\nu,2},\dots,z_{\nu,\phi(\nu)}$
for the zeros of $K_\nu$ and it is known that $\RE(z_{\nu,j})<0$
for each $j$. Since $K_\nu$ is one of the fundamental solutions of
the second order equation (2.2), it is easy to see that
all zeros of $K_\nu$ are of multiplicity one
by the uniqueness of the solution of ordinary differential equations.
This yields that all zeros of $K_\nu$ are distinct.
For details, see \cite{18, pp.511--513}. This article deals with
the case of $\nu=d/2-1$. For convenience, we use $z_j^{(d)}$ and $\phi_d$
instead of $z_{d/2-1,j}$ and $\phi(d/2-1)$, respectively.

When $d$ is odd and not less than five,
Theorem 1.1 in \cite{6} shows that, for $t>0$
$$
L(t)=S_r^{(d)}
\biggl[\frac{(d-2)t}2+\frac{r^2}{d-4}
-\frac{\sqrt2 r^3}{\sqrt{\pi t}}\sum_{j=1}^{\phi_d}
\frac1{(z_j^{(d)})^2}\int_0^\infty 
e^{-\frac{r^2x^2}{2t}+z_j^{(d)} x}dx\biggr].
\tag2.3
$$
We have used $S_r^{(d)}=S_{d-1}r^{d-2}$ for simplicity.
When $d$ is even, Theorem 4.1 in \cite{9} shows that, for $t>0$,
if $d=2$,
$$
L(t)=2\pi r\biggl[\sqrt{\frac{2t}\pi}+\frac{\sqrt 2 r^2}{\sqrt{\pi t}}
\int_0^\infty \!\!\! \int_0^\infty \frac{xy-1+e^{-xy}}{y^3 G^{(2)}(y)}
e^{-\frac{r^2x^2}{2t}}dxdy\biggr],
$$
if $d=4$,
$$
L(t)=2\pi^2r^2\biggl[t+\frac{\sqrt 2 r^3}{\sqrt{\pi t}}
\int_0^\infty \!\!\! \int_0^\infty \frac{1-e^{-xy}}{y^3 G^{(4)}(y)}
e^{-\frac{r^2x^2}{2t}}dxdy\biggr],
\tag2.4
$$
and, if $d\geqq6$,
$$
\split
L(t)=S_r^{(d)}\biggl[&\frac{(d-2)t}2+\frac{r^2}{d-4}\\
&-\frac{\sqrt 2 r^3}{\sqrt{\pi t}}\sum_{j=1}^{\phi_d}
\frac1{(z_j^{(d)})^2}\int_0^\infty e^{-\frac{r^2x^2}{2t}+z_j^{(d)}x}dx\\
&+\frac{(-1)^{d/2-1}\sqrt 2 r^3}{\sqrt{\pi t}}
\int_0^\infty \!\!\! \int_0^\infty \!\!\frac{e^{-xy}}{y^3 G^{(d)}(y)}
e^{-\frac{r^2x^2}{2t}}dxdy\biggr].
\endsplit
\tag2.5
$$
The notation $G^{(d)}(x)$ has been used to denote
$K_{d/2-1}(x)^2+\pi^2 I_{d/2-1}(x)^2$ and, if not confuse,
we use $G(x)$ instead of $G^{(d)}(x)$ .

The remainder of this section is devoted to the asymptotic expansion
of $L(t;A)$ for large $t$. Let $M\geqq1$ be a given integer.
In the two dimensional case, it has been already proved
$$
L(t;A)=\frac{2\pi t}{\log t}\biggl[
1+\frac{c_1^{(2)}}{\log t}+\frac{c_2^{(2)}}{(\log t)^2}+
\dots+\frac{c_M^{(2)}}{(\log t)^M}\biggr]
+O\mm{\frac t{(\log t)^{M+2}}}
$$
for some sequence $\{c_n^{(2)}\}_{n=1}^\infty$ of real numbers
(cf. \cite{15}). In the higher dimensional case, it has been shown
$$
L(t;A)=
\cases
c_1^{(4)}t+c_2^{(4)}\log t+c_3^{(4)}
+c_4^{(4)}\dfrac{\log t}t+o\mm{\dfrac{\log t}t}
&\quad\text{if $d=4$,}\\
c_1^{(d)}t+c_2^{(d)}+c_3^{(d)}t^{2-d/2}
+O\mm{t^{1-d/2}}&\quad\text{if $d\geqq5$}
\endcases
\tag2.6
$$
and been given each constant $c_j^{(d)}$ explicitly
(cf. \cite{14}). In the case of $A=D$, we can improve (2.6).
When $d$ is odd and not less than five, Theorem 4.1 in \cite{6} shows
$$
\split
L(t)=
S_r^{(d)}\biggl[&\frac{(d-2)t}2+\frac{r^2}{d-4}\\
&+\sum_{n=(d-5)/2}^M
\sum_{j=1}^{\phi_d}\frac{\lambda_n r^{2n+3}}{(z_j^{(d)})^{2n+3}}
\frac1{t^{n+1/2}}\biggr]+O\biggl[\frac1{t^{M+3/2}}\biggr],
\endsplit
\tag2.7
$$
where
$$
\lambda_n=\frac{(-1)^n\varGamma(2n+1)}{2^{n-1/2}\sqrt\pi \varGamma(n+1)}
$$
and $\varGamma$ is the gamma function.

If $d$ is even and $d\geqq6$, the third term of (2.5)
is the same as that of (2.3).
This means that the expansion of the third term of (2.5) can be derived
by the method used to obtain (2.7).
Hence we need to consider the double integral
in (2.5) and have the following result whose proof is
deferred to Section~4.
For convenience, we put $N=d/2-1$. Note that $N$ is an integer
and not less than $2$.

\medskip

\proclaim{Theorem 2.1}
Let $M\geqq1$ be a given integer. If $d\geqq6$ and $d$ is even, we have
$$
\split
L(t)=&S_r^{(d)}\biggl[\frac{d-2}2t+\frac{r^2}{d-4}
+\sum_{m=0}^{M-1}\sum_{n=0}^{N-1}\frac{r^{2N(m+1)+2n}}{t^{N(m+1)+n-1}}
\sum_{k=0}^m\rho_{m,n,k}^{(1)}(\log t)^k\\
&\hphantom{S_r^{(d)}\biggl[}+\sum_{m=0}^{M-1}
\sum_{n=0}^{N(M-m)-1}\frac{r^{2N(m+1)+2n+1}}{t^{N(m+1)+n-1/2}}
\sum_{k=0}^m\rho_{m,n,k}^{(2)}(\log t)^k\\
&\hphantom{S_r^{(d)}\biggl[}+\frac{r^{2N(M+1)}}{t^{N(M+1)-1}}
\sum_{k=0}^M\rho_k^{(3)}(\log t)^k
+\frac{r^{2N(M+1)+1}}{t^{N(M+1)-1/2}}
\sum_{k=0}^M\rho_k^{(4)}(\log t)^k\biggr]\\
&+O\mm{\frac1{t^{N(M+1)-1/2}}}
\endsplit
$$
for sequences $\{\rho_{n,m,k}^{(1)}\}$,
$\{\rho_{n,m,k}^{(2)}\}$, $\{\rho_k^{(3)}\}$,
$\{\rho_k^{(4)}\}$ of real numbers.
\endproclaim

\medskip

If $d=4$, the calculation of (2.4) is similar. However, since $G(x)$
is asymptotically equal to $1/x^2$ as $x\downarrow0$,
which will be discussed in the next section, we must treat
the double integral in (2.4) with much care.
The result in the four dimensional case is the following,
which will be proved in Section 5.

\medskip

\proclaim{Theorem 2.2}
Let $M\geqq1$ be a given integer. If $d=4$, we have
$$
\split
L(t)=2\pi^2r^2\biggl[t+\frac12 r^2\log t&+r^2\sigma^{(1)}
+\sum_{n=1}^M\frac{r^{2n+2}}{t^n}\sum_{k=0}^n\sigma_{n,k}^{(2)}(\log t)^k\\
&+\sum_{n=1}^M\frac{r^{2n+3}}{t^{n+1/2}}\sum_{k=0}^n
\sigma_{n,k}^{(3)}(\log t)^k
\biggr]+O\mm{\frac{(\log t)^{M+1}}{t^{M+1}}}
\endsplit
$$
for a constant $\sigma^{(1)}$ and suitable
sequences $\{\sigma_{n,k}^{(2)}\}$, $\{\sigma_{n,k}^{(3)}\}$.
\endproclaim

%333333333333333333333333333333333333333

\head
3. Modified Bessel functions
\endhead

\noindent
This section is devoted to the asymptotic expansion of $G(x)$
for small $x$ in the case when $d$ is even and larger than or equal
to $4$, which will be used to compute the double integral in (2.5).
Recall $N=d/2-1$. We have that $N$ is a positive integer and that
$$
G(x)=K_N(x)^2+\pi^2 I_N(x)^2.
$$
With the help of asymptotic behaviors of $K_N$ and $I_N$,
it can be easily shown
$$
G(x)=
\cases
\dfrac1{\kappa_N x^{2N}}\{1+o[1]\}\quad&\text{as $x\downarrow0$,}\\
\dfrac{\pi e^{2x}}{2x}\{1+o[1]\}\quad&\text{as $x\to\infty$,}
\endcases
\tag3.1
$$
where the notation $\kappa_N=1/4^{N-1}\varGamma(N)^2$ (cf. \cite{9}).
However (3.1) for small $x$ is not sufficient for our purpose
and thus we need to improve it in the following way.

\medskip

\proclaim{Lemma 3.1}
Let $L\geqq1$ be a given integer. We have that, as $x\downarrow0$,
$$
\split
\kappa_Nx^{2N}G(x)=1&-x^2\sum_{n=0}^L\alpha_n x^{2n}
-x^{2N}\log x\sum_{n=0}^L\beta_n x^{2n}\\
&-x^{4N}(\log x)^2\sum_{n=0}^L\gamma_n x^{2n}+O\mm{x^{2(L+2)}\log x}
\endsplit
$$
for some sequences $\{\alpha_n\}_{n=0}^\infty$,
$\{\beta_n\}_{n=0}^\infty$, $\{\gamma_n\}_{n=0}^\infty$
of real numbers.
\endproclaim
\demo{Proof}
It is sufficient to prove this lemma for $L\geqq 2N-1$.
For $n\geqq0$ let
$$
\split
&\xi_n=
\cases
\dfrac{(-1)^n\varGamma(N-n)}
{2^{2n-N+1}\varGamma(n+1)}&\text{if $0\leqq n\leqq N-1$,}\\
\dfrac{(-1)^N\{\log4+\psi(n-N+1)+\psi(n+1)\}}
{2^{2n-N+1}\varGamma(n-N+1)\varGamma(n+1)}\quad
&\text{if $n\geqq N$,}
\endcases\\
&\eta_n=
\frac{(-1)^{N+1}}{2^{2n+N}\varGamma(n+1)\varGamma(n+N+1)},
\endsplit
$$
where $\psi$ is the first derivative of $\log \varGamma$.
The formula for $K_N$ yields
$$
x^NK_N(x)=
\sum_{n=0}^{N-1} \xi_n x^{2n}
+x^{2N}\sum_{n=0}^\infty \xi_{n+N} x^{2n}
+x^{2N}\log x\sum_{n=0}^\infty \eta_n x^{2n}
$$
(see \cite{18, p.80}). Hence we have that, as $x\downarrow0$,
$$
x^NK_N(x)=\sum_{n=0}^{L+1}\xi_n x^{2n}
+x^{2N}\log x\sum_{n=0}^{L+1}\eta_nx^{2n}
+O\mm{x^{2(L+2)}}.
$$
This implies that $x^{2N}K_N(x)^2$ is
$$
\split
&\sum_{k=0}^{L+1}\sum_{n=0}^{L+1}\xi_k\xi_n x^{2(k+n)}
+2x^{2N}\log x\sum_{k=0}^{L+1}\sum_{n=0}^{L+1}\xi_k\eta_n x^{2(k+n)}\\
&\hphantom{-}
+x^{4N}(\log x)^2\sum_{k=0}^{L+1}\sum_{n=0}^{L+1}\eta_k\eta_n x^{2(k+n)}
+O\mm{x^{2(L+2)}}.
\endsplit
\tag3.2
$$
It is easy to see that the first term of (3.2) is
$$
\sum_{k=0}^{L+1}\sum_{n=k}^{L+1}\xi_k\xi_{n-k}x^{2n}+O\mm{x^{2(L+2)}}
=\sum_{n=0}^{L+1}\biggl(\sum_{k=0}^n\xi_k\xi_{n-k}\biggr)x^{2n}
+O\mm{x^{2(L+2)}},
$$
which is equal to
$$
\xi_0^2+x^2\sum_{n=0}^L\biggl(\sum_{k=0}^{n+1}\xi_k\xi_{n-k+1}\biggr)x^{2n}
+O\mm{x^{2(L+2)}}.
$$
Similarly we have that the second and the third terms of (3.2) are
$$
x^{2N}\log x\sum_{n=0}^{L+1}\biggl(2\sum_{k=0}^n\xi_k\eta_{n-k}\biggr)x^{2n}
+O\mm{x^{2(L+N+2)}\log x}
$$
and
$$
x^{4N}(\log x)^2\sum_{n=0}^{L+1}\biggl(\sum_{k=0}^n\eta_k\eta_{n-k}\biggr)x^{2n}
+O\mm{x^{2(L+N+2)}(\log x)^2},
$$
respectively. Note that $\xi_0^2=1/\kappa_N$ and then it follows that
$$
\split
x^{2N}K_N(x)^2=
\frac1{\kappa_N}&+x^2\sum_{n=0}^L
\biggl(\sum_{k=0}^{n+1}\xi_k\xi_{n-k+1}\biggr)x^{2n}\\
&+x^{2N}\log x\sum_{n=0}^L
\biggl(2\sum_{k=0}^n\xi_k\eta_{n-k}\biggr)x^{2n}\\
&+x^{4N}(\log x)^2\sum_{n=0}^L
\biggl(\sum_{k=0}^n\eta_k\eta_{n-k}\biggr)x^{2n}\\
&+O\mm{x^{2(L+2)}+x^{2(L+N+1)}\log x}.
\endsplit
$$

On the other hand, the formula for $I_N$ yields
$$
\pi x^NI_N(x)=x^{2N}\sum_{n=0}^{L-2N+1} \zeta_n x^{2n}
+O\mm{x^{2(L-N+2)}},
\tag3.3
$$
where
$$
\zeta_n=\frac\pi{2^{2n+N}\varGamma(n+1)\varGamma(n+N+1)}
$$
for $n\geqq0$ (see \cite{18, p.77}). Since $2(L-N+2)\geqq L+2$,
we deduce from (3.3) that
$$
\split
\pi^2x^{2N}I_N(x)^2
&=x^{4N}\sum_{k=0}^{L-2N+1}\sum_{n=k}^{L-2N+1}\zeta_k\zeta_{n-k}x^{2n}
+O\mm{x^{2(L+2)}}\\
&=x^{4N}\sum_{n=0}^{L-2N+1}
\biggl(\sum_{k=0}^n\zeta_k\zeta_{n-k}\biggr)x^{2n}
+O\mm{x^{2(L+2)}},
\endsplit
$$
which is equal to
$$
x^2\sum_{n=2N-1}^L\biggl(\sum_{k=0}^{n-2N+1}
\zeta_k\zeta_{n-2N-k+1}\biggr)x^{2n}
+O\mm{x^{2(L+2)}}.
$$
For an integer $n\geqq0$ let
$$
\split
&\alpha_n=
\cases
\dis-\kappa_N\sum_{k=0}^{n+1}\xi_k\xi_{n-k+1}
&\text{if $0\leqq n\leqq 2N-2$,}\\
\dis-\kappa_N\biggl(\sum_{k=0}^{n+1}\xi_k\xi_{n-k+1}
+\sum_{k=0}^{n-2N+1}\zeta_k\zeta_{n-2N-k+1}\biggr)\quad
&\text{if $n\geqq 2N-1$,}
\endcases\\
&\beta_n=-2\kappa_N\sum_{k=0}^n\xi_k\eta_{n-k},\quad
\gamma_n=-\kappa_N\sum_{k=0}^n\eta_k\eta_{n-k}.\\
\endsplit
$$
Hence we have
$$
\split
\kappa_Nx^{2N}G(x)=1&-x^2\sum_{n=0}^L\alpha_n x^{2n}
-x^{2N}\log x\sum_{n=0}^L\beta_n x^{2n}\\
&-x^{4N}(\log x)^2\sum_{n=0}^L\gamma_n x^{2n}
+\cases
O\mm{x^{2(L+2)}\log x}\quad&\text{if $N=1$,}\\
O\mm{x^{2(L+2)}}&\text{if $N\geqq2$.}
\endcases
\endsplit
$$
This completes the proof of this lemma.
\qed
\enddemo

\medskip

Let $\{p_{n,0}\}_{n=0}^\infty$ and $\{p_{n,1}\}_{n=0}^\infty$
be sequences of real numbers defined by
$$
\cases
p_{0,0}=1,\\
\dis p_{n,0}=\sum_{m=0}^{n-1}\alpha_{n-m-1}p_{m,0}\quad(n\geqq0)
\endcases
$$
and
$$
\cases
p_{0,1}=\beta_0,\\
\dis p_{n,1}=\sum_{m=0}^{n-1}\alpha_{n-m-1}p_{m,1}
+\sum_{m=0}^n\beta_{n-m}p_{m,0}\quad(n\geqq0).
\endcases
$$
Under the assumption that $\{p_{n,h}\}_{n=0}^\infty$
have been determined for $0\leqq h\leqq k-1$,
we define a sequence $\{p_{n,k}\}_{n=0}^\infty$ of real numbers by
$$
\cases
p_{0,k}=\beta_0p_{0,k-1}+\gamma_0p_{0,k-2},\\
\dis p_{n,k}=\sum_{m=0}^{n-1}\alpha_{n-m-1}p_{m,k}
+\sum_{m=0}^n\beta_{n-m}p_{m,k-1}
+\sum_{m=0}^n\gamma_{n-m}p_{m,k-2}\quad(n\geqq0).
\endcases
$$
For integers $n$, $j$, $k$ with $0\leqq n\leqq N-1$, $0\leqq k\leqq j$
let
$$
p_{n,k}^j=p_{n+(j-k)N,k}.
\tag3.4
$$
Using the sequence $\{p_{n,k}^j\}$,
we can obtain the asymptotic expansion of $1/G(x)$ as $x\downarrow0$.

\medskip

\proclaim{Proposition 3.2}
For a given integer $M\geqq0$ we have that, as $x\downarrow0$,
$$
\frac1{\kappa_Nx^{2N}G(x)}=\sum_{j=0}^M\sum_{n=0}^{N-1}x^{2n+2Nj}
\sum_{k=0}^j p_{n,k}^j(\log x)^k+O\mm{x^{2N(M+1)}(\log x)^{M+1}}.
$$
\endproclaim
\demo{Proof}
It is sufficient to see the theorem for $M\geqq3$. Moreover
we consider only the case of $N\geqq2$ since the case of $N=1$ can be
shown in the similar way. Hence we give a proof for $M\geqq3$
when $N\geqq2$.

Since $\kappa_Nx^{2N}G(x)=1+o[1]$ (see (3.1)),
it is sufficient to see that
$$
\kappa_Nx^{2N}G(x)\sum_{j=0}^M\sum_{n=0}^{N-1}x^{2n+2Nj}
\sum_{k=0}^j p_{n,k}^j (\log x)^k
\tag3.5
$$
is $1+O[x^{2N(M+1)} (\log x)^{M+1}]$ as $x\downarrow0$.
Applying Lemma 3.1 for $L=N(M+1)-1$, we obtain
$$
\split
\kappa_Nx^{2N}G(x)=1&-x^2\sum_{n=0}^{N(M+1)-1}\alpha_n x^{2n}
-x^{2N}\log x\sum_{n=0}^{NM}\beta_n x^{2n}\\
&-x^{4N}(\log x)^2\sum_{n=0}^{N(M-1)}\gamma_n x^{2n}
+O\mm{x^{2N(M+1)+2}(\log x)^2}.
\endsplit
$$
This formula implies that the main part of (3.5) is
$$
\align
&\sum_{j=0}^M \sum_{n=0}^{N-1}x^{2n+2Nj}
\sum_{k=0}^j p_{n,k}^j (\log x)^k\tag3.6\\
&\hphantom{-}-\sum_{j=0}^M \sum_{m=0}^{N-1} \sum_{n=0}^{N(M+1)-1}
\alpha_n x^{2(m+n+1)+2Nj}\sum_{k=0}^j p_{m,k}^j (\log x)^k\tag3.7\\
&\hphantom{-}-\sum_{j=0}^M \sum_{m=0}^{N-1} \sum_{n=0}^{NM}
\beta_n x^{2(m+n)+2N(j+1)}\sum_{k=0}^j p_{m,k}^j (\log x)^{k+1}\tag3.8\\
&\hphantom{-}-\sum_{j=0}^M \sum_{m=0}^{N-1} \sum_{n=0}^{N(M-1)}
\gamma_n x^{2(m+n)+2N(j+2)}\sum_{k=0}^j p_{m,k}^j (\log x)^{k+2}\tag3.9\\
\endalign
$$
and the remaining part of (3.5) is $O[x^{2N(M+1)}(\log x)^{M+1}]$.

\pagebreak

It follows from $p_{0,0}^0=1$ that (3.6) is equal to
$$
1+\sum_{n=1}^{N-1}x^{2n}p_{n,0}^0+\sum_{j=1}^M\sum_{n=0}^{N-1}x^{2n+2Nj}
\sum_{k=0}^j p_{n,k}^j (\log x)^k.
$$
For the calculation of (3.7) we first change the index of
the summation from $n$ to $l$ given by $l=n+m+1$ and next write $n$ again instead of $l$.
Then we have that (3.7) is
$$
-\sum_{j=0}^M \sum_{m=0}^{N-1} \sum_{n=m+1}^{N(M+1)+m}\alpha_{n-m-1}
x^{2n+2Nj}\sum_{k=0}^j p_{m,k}^j (\log x)^k.
$$
Since the contribution of the summation on $n$ over $[N(M+1), N(M+1)+m]$
is of order $x^{2N(M+1)}$, the leading term of (3.7) is equal to
$$
-\sum_{j=0}^M\sum_{m=0}^{N-1}\sum_{n=m+1}^{N(M+1)-1}\alpha_{n-m-1}x^{2n+2Nj}
\sum_{k=0}^j p_{m,k}^j (\log x)^k.
$$
We divide the summation on $n$ into two parts.
One is the sum over $[m+1, N]$ and the other is
the sum over $[N+1, N(M+1)-1]$.
Therefore, changing the order of the summations on $n$ and $m$,
we derive that the first part is
$$
-\sum_{j=0}^M\sum_{n=1}^N x^{2n+2Nj}\sum_{k=0}^j
\biggl(\sum_{m=0}^{n-1}\alpha_{n-m-1}p_{m,k}^j\biggr)(\log x)^k
\tag3.10
$$
and that the second one is
$$
-\sum_{j=0}^M\sum_{n=N+1}^{N(M+1)-1} x^{2n+2Nj}\sum_{k=0}^j
\biggl(\sum_{m=0}^{N-1}\alpha_{n-m-1}p_{m,k}^j\biggr)(\log x)^k.
\tag3.11
$$
Note that the case of $n=N$ in (3.10) coincides with that in (3.11).
Hence (3.7) is
$$
\split
-&\sum_{j=0}^M\sum_{n=1}^{N-1} x^{2n+2Nj}\sum_{k=0}^j
\biggl(\sum_{m=0}^{n-1}\alpha_{n-m-1}p_{m,k}^j\biggr)(\log x)^k\\
&-\sum_{j=0}^M\sum_{n=N}^{N(M+1)-1} x^{2n+2Nj}\sum_{k=0}^j
\biggl(\sum_{m=0}^{N-1}\alpha_{n-m-1}p_{m,k}^j\biggr)(\log x)^k\\
&+O\mm{x^{2N(M+1)}}.
\endsplit
\tag3.12
$$
We must still calculate the second term of (3.12).
Dividing the summation on $n$ into $M$ parts such that
the $i$-th part is the sum over $[iN, (i+1)N-1]$, we obtain that
the second term of (3.12) is
$$
-\sum_{j=0}^M\sum_{i=1}^M\sum_{n=iN}^{(i+1)N-1} x^{2n+2Nj}\sum_{k=0}^j
\biggl(\sum_{m=0}^{N-1}\alpha_{n-m-1}p_{m,k}^j\biggr)(\log x)^k,
$$
which coincides with
$$
\sum_{i=1}^M\sum_{j=0}^M\sum_{n=0}^{N-1} x^{2n+2N(i+j)}\sum_{k=0}^j
\biggl(\sum_{m=0}^{N-1}\alpha_{n+iN-m-1}p_{m,k}^j\biggr)(\log x)^k.
\tag3.13
$$
Here we first have changed the index of the sum from $n$ to $l=n+iN$
and next have written $n$ again instead of $l$.
Since the contribution of the summation on $j$ over $[M-i,M]$
in (3.13) is of order $x^{2N(M+1)}(\log x)^M$,
the leading term of (3.13) is
$$
-\sum_{i=1}^M\sum_{j=i}^M\sum_{n=0}^{N-1} x^{2n+2Nj}\sum_{k=0}^{j-i}
\biggl(\sum_{m=0}^{N-1}\alpha_{n+iN-m-1}p_{m,k}^{j-i}\biggr)(\log x)^k.
$$
A change of the order of the summations on $i$, $j$, $k$
yields that (3.12) and also (3.7) are equal to
$$
\split
-&\sum_{j=0}^M\sum_{n=1}^{N-1} x^{2n+2Nj}\sum_{k=0}^j
\biggl(\sum_{m=0}^{n-1}\alpha_{n-m-1}p_{m,k}^j\biggr)(\log x)^k\\
&-\sum_{j=1}^M\sum_{n=0}^{N-1}x^{2n+2Nj}\sum_{k=0}^{j-1}
\biggl(\sum_{i=1}^{j-k}\sum_{m=0}^{N-1}\alpha_{n+iN-m-1}
p_{m,k}^{j-i}\biggr)(\log x)^k\\
&+O\mm{x^{2N(M+1)}(\log x)^M}.
\endsplit
$$

The calculation of (3.8) is similar but easier than that of (3.7).
The contribution of the sum on $n$ over $[NM-m,NM]$
is of order $x^{2N(M+1)}(\log x)^{M+1}$. Then the main term of
(3.8) is
$$
-\sum_{j=1}^M\sum_{m=0}^{N-1}\sum_{n=m}^{NM-1}
\beta_{n-m}x^{2n+2Nj}\sum_{k=1}^jp_{m,k-1}^j(\log x)^k.
\tag3.14
$$
We divide the summation on $n$ into two parts such that
the sums over $[m, N-1]$ and $[N, NM-1]$.
Hence, changing the order of
the summations on $n$ and $m$, we have that (3.14) is
$$
\split
-&\sum_{j=1}^M\sum_{n=0}^{N-1} x^{2n+2Nj}\sum_{k=1}^j
\biggl(\sum_{m=0}^n \beta_{n-m}p_{m,k-1}^{j-1}\biggr)(\log x)^k\\
&-\sum_{j=1}^M\sum_{n=N}^{NM-1} x^{2n+2Nj}\sum_{k=1}^j
\biggl(\sum_{m=0}^{N-1}\beta_{n-m}p_{m,k-1}^{j-1}\biggr)(\log x)^k.\\
\endsplit
\tag3.15
$$
Similarly to the second term of (3.12), we can derive that the second
term of (3.15) is the sum of
$$
-\sum_{j=2}^M\sum_{n=0}^{N-1}x^{2n+2Nj}\sum_{k=1}^{j-1}
\biggl(\sum_{i=1}^{j-k}\sum_{m=0}^{N-1}
\beta_{n+iN-m}p_{m,k-1}^{j-i-1}\biggr) (\log x)^k
$$
and the part which is of order $x^{2N(M+1)}(\log x)^M$.
The same calculation as (3.8) shows that the leading term of (3.9) is
$$
\split
-&\sum_{j=2}^M\sum_{n=0}^{N-1} x^{2n+2Nj}\sum_{k=2}^j
\biggl(\sum_{m=0}^n \gamma_{n-m}p_{m,k-2}^{j-2}\biggr)(\log x)^k\\
&-\sum_{j=3}^M\sum_{n=0}^{N-1}x^{2n+2Nj}\sum_{k=2}^{j-1}\biggl(
\sum_{i=1}^{j-k}\sum_{m=0}^{N-1}\gamma_{n+iN-m}p_{m,k-2}^{j-i-2}\biggr)
(\log x)^k\\
\endsplit
$$
and its error term is $O[x^{2N(M+1)}(\log x)^{M+1}]$.
Details are left to the reader.

Combining these calculations, we conclude that the main part of (3.5) is
$$
\split
1&+\sum_{n=1}^{N-1}x^{2n}p_{n,0}^0
+\sum_{j=1}^M\sum_{n=0}^{N-1}x^{2n+2Nj}
\sum_{k=0}^j p_{n,k}^j (\log x)^k\\
&-\sum_{j=0}^M\sum_{n=1}^{N-1} x^{2n+2Nj}\sum_{k=0}^j
\sum_{m=0}^{n-1}\alpha_{n-m-1}p_{m,k}^j(\log x)^k\\
&-\sum_{j=1}^M\sum_{n=0}^{N-1}x^{2n+2Nj}\sum_{k=0}^{j-1}
\sum_{i=1}^{j-k}\sum_{m=0}^{N-1}\alpha_{n+iN-m-1}
p_{m,k}^{j-i}(\log x)^k\\\allowdisplaybreak
&-\sum_{j=1}^M\sum_{n=0}^{N-1} x^{2n+2Nj}\sum_{k=1}^j
\sum_{m=0}^n \beta_{n-m}p_{m,k-1}^{j-1}(\log x)^k\\
&-\sum_{j=2}^M\sum_{n=0}^{N-1}x^{2n+2Nj}\sum_{k=1}^{j-1}
\sum_{i=1}^{j-k}\sum_{m=0}^{N-1}\beta_{n+iN-m}p_{m,k-1}^{j-i-1}
(\log x)^k\\
&-\sum_{j=2}^M\sum_{n=0}^{N-1} x^{2n+2Nj}\sum_{k=2}^j
\sum_{m=0}^n \gamma_{n-m}p_{m,k-2}^{j-2}(\log x)^k\\
&-\sum_{j=3}^M\sum_{n=0}^{N-1}x^{2n+2Nj}\sum_{k=2}^{j-1}
\sum_{i=1}^{j-k}\sum_{m=0}^{N-1}\gamma_{n+iN-m}p_{m,k-2}^{j-i-2}
(\log x)^k\\
\endsplit
\tag3.16
$$
and the remaining part of (3.5) is of order $x^{2N(M+1)}(\log x)^{M+1}$.

We now concentrate on showing that all coefficients in (3.16) varnish
except for the constant. It is sufficient to consider
only the coefficient of $x^{2n+2Nj}(\log x)^k$
for $1\leqq n\leqq N-1$, $j\geqq3$ and $2\leqq k\leqq j-1$
since the way of calculations of other terms are similar.
We deduce from (3.4) that the coefficient of $x^{2n+2Nj}(\log x)^k$
in (3.16) is equal to
$$
\split
&p_{n+(j-k)N,k}-\sum_{m=0}^{n-1}\alpha_{n-m-1}p_{m+(j-k)N,k}
-\sum_{i=1}^{j-k}\sum_{m=0}^{N-1}\alpha_{n+iN-m-1}p_{m+(j-k-i)N,k}\\
&\hphantom{=}-\sum_{m=0}^n\beta_{n-m}p_{m+(j-k)N,k-1}
-\sum_{i=1}^{j-k}\sum_{m=0}^{N-1}\beta_{n+iN-m}p_{m+(j-k-i)N,k-1}\\
&\hphantom{=}-\sum_{m=0}^n\gamma_{n-m}p_{m+(j-k)N,k-2}
-\sum_{i=1}^{j-k}\sum_{m=0}^{N-1}\gamma_{n+iN-m}p_{m+(j-k-i)N,k-2}.\\
\endsplit
$$
It is easy to see that the second term is
$$
\sum_{m=(j-k)N}^{n+(j-k)N-1}\alpha_{n+(j-k)N-m-1}p_{m,k}.
$$
Similarly we derive that the fourth and the sixth terms are
$$
\sum_{m=(j-k)N}^{n+(j-k)N}\beta_{n+(j-k)N-m}p_{m,k-1},\quad
\sum_{m=(j-k)N}^{n+(j-k)N}\gamma_{n+(j-k)N-m}p_{m,k-2},
$$
respectively. Moreover the third term is
$$
\sum_{i=1}^{j-k}\sum_{m=(j-k-i)N}^{(j-k-i+1)N-1}
\alpha_{n+(j-k)N-m-1}p_{m,k}
=\sum_{m=0}^{(j-k)N-1}\alpha_{n+(j-k)N-m-1}p_{m,k}.
$$
The same calculation shows that the fifth and the last terms are equal to
$$
\sum_{m=0}^{(j-k)N-1}\beta_{n+(j-k)N-m}p_{m,k-1},\quad
\sum_{m=0}^{(j-k)N-1}\gamma_{n+(j-k)N-m}p_{m,k-2},
$$
respectively. Therefore we have that the coefficient of
$x^{2n+2Nj}(\log x)^k$ in (3.16) coincides with
$$
\split
&p_{n+(j-k)N,k}-\sum_{m=0}^{n+(j-k)N-1}\alpha_{n+(j-k)N-m-1}p_{m+(j-k)N,k}\\
&\hphantom{=}-\sum_{m=0}^{n+(j-k)N}\beta_{n+(j-k)N-m}p_{m+(j-k)N,k-1}\\
&\hphantom{=}-\sum_{m=0}^{n+(j-k)N}\gamma_{n+(j-k)N-m}p_{m+(j-k)N,k-2}.\\
\endsplit
$$
The definition of $p_{n,k}$ yields that it varnishes.
\qed
\enddemo

%44444444444444444444444444444444444444444

\head
4. The six or more dimensional cases
\endhead

\noindent
Our goal in this section is to show Theorem 2.1 with the help of (2.5).
We assume that $d$ is even and not less than 6. Recall that we put
$N=d/2-1$ in Section 2 and thus $N$ is an integer
which is larger than or equal to 2.
Throughout this section, we use $C_1,\dots,C_7$ for 
positive constants independent of the variable.

We first consider the fourth term of the right hand side of (2.5).
For $t>0$ let
$$
S(t)=\frac{\sqrt2 r^3}{\sqrt{\pi t}}
\int_0^\infty \!\!\! \int_0^\infty \frac{e^{-xy}}{y^3 G(y)}
e^{-\frac{r^2x^2}{2t}}dxdy.
$$
We set $T(t)=S(r^2t^2/2)$ and thus obtain that
$$
T(t)=\frac{2r^2}{\sqrt\pi\,t}\int_0^\infty\!\!\! \int_0^\infty
\frac{e^{-xy}}{y^3 G(y)} e^{-\frac{x^2}{t^2}}dxdy.
$$
A change of the variable from $x$ to $u=x/t$ shows
that $T(t)$ is the sum of
$$
\align
&\frac{2r^2}{\sqrt\pi}\int_1^\infty dy\int_0^\infty
\frac1{y^3G(y)}e^{-tuy} e^{-u^2} du,\tag4.1\\
&\frac{2r^2}{\sqrt\pi}\int_0^1 dy\int_0^\infty
\frac1{y^3G(y)}e^{-tuy} e^{-u^2} du.\tag4.2
\endalign
$$

We first try to calculate (4.1). For simplicity, we put
$$
a_m=\frac{(-1)^m}{\varGamma(m+1)},
\quad
b_m=\frac{(-1)^m\varGamma(2m+1)}{\varGamma(m+1)}
$$
for an integer $m\geqq0$. For $x>0$ and an integer $n\geqq0$ let
$$
R_n(x)=e^{-x}-\sum_{k=0}^n a_k x^k.
$$
It is easy to see that, for $n\geqq0$,
$$
\align
&|R_n(x)|\leqq A_n|x|^{n+1},\quad x>0,\tag4.3\\
&|R_n(x)|\leqq B_n|x|^n,\quad x\geqq1\tag4.4
\endalign
$$
for some suitable constants $A_n$ and $B_n$.
Applying the well-known formula
$$
\int_0^1 x^{a-1}\biggl(\log\frac1x\biggr)^pdx=\int_0^\infty e^{-ax} x^p dx
=\frac{\varGamma(p+1)}{a^{p+1}}
\tag4.5
$$
for $p>-1$ and $a>0$, we have by (4.3) that
$$
\split
&\biggl|\int_1^\infty dy \int_0^\infty \frac1{y^3G(y)}e^{-tuy}
R_{N(M+1)-2}(u^2) du\biggr|\\
&\leqq \frac{C_1}{t^{2N(M+1)-1}}
\int_1^\infty \frac{dy}{y^{2N(M+1)+2}G(y)}.
\endsplit
$$
Note that (3.1) yields that $1/y^p G(y)$ is integrable on $[1,\infty)$
for any $p\in\Bbb R$. Hence it follows that (4.1) is the sum of
$$
\frac{2r^2}{\sqrt\pi}\sum_{n=0}^{N(M+1)-2} a_n
\int_1^\infty dy\int_0^\infty \frac1{y^3G(y)}e^{-tuy} u^{2n} du
\tag4.6
$$
and the part which is of order $1/t^{2N(M+1)-1}$. Applying (4.5) again,
we have that (4.6) is
$$
\frac{2r^2}{\sqrt\pi}\sum_{n=0}^{N(M+1)-2} b_n
\int_1^\infty \frac{dy}{y^{2n+4}G(y)}\frac1{t^{2n+1}}
\tag4.7
$$
and thus that the main term of (4.1) is (4.7) and the remaining term
of (4.1) is of order $1/t^{2N(M+1)-1}$.

We next calculate (4.2), which is equal to
$$
\align
&\frac{2r^2}{\sqrt\pi}\sum_{n=0}^{N-2} a_n
\int_0^1 dy\int_0^\infty \frac1{y^3G(y)}e^{-tuy}u^{2n}du
\tag4.8\\
&\hphantom{-}+\frac{2r^2}{\sqrt\pi} \int_0^1 dy \int_0^\infty
\frac1{y^3G(y)}e^{-tuy}R_{N-2}(u^2)du.
\tag4.9
\endalign
$$
We remark that $1/y^pG(y)$ is integrable on $(0,1)$
for $p<2N+1$, which can be seen from (3.1). The same calculation
as (4.6) yields that (4.8) is
$$
\frac{2r^2}{\sqrt\pi}\sum_{n=0}^{N-2} b_n
\int_0^1 \frac{dy}{y^{2n+4}G(y)}\frac1{t^{2n+1}}.
$$
For the calculation of (4.9)
we need to consider $R_N(u^2)/G(y)$. For $0\leqq h\leqq N-1$ let
$$
Q_h^0(x)=\frac1{G(x)}-\kappa_N x^{2N}\sum_{n=0}^h x^{2n}p_n^0(x)
$$
and for $0\leqq h\leqq N-1$, $m\geqq1$ let
$$
\split
Q_h^m(x)=\frac1{G(x)}
&-\kappa_N x^{2N}\sum_{j=0}^{m-1}\sum_{n=0}^{N-1} x^{2n+2Nj}p_n^j(x)\\
&-\kappa_N x^{2N}\sum_{n=0}^h x^{2n+2Nm}p_n^m(x).
\endsplit
$$
Here we have used the notation $p_n^j(x)$ defined by
$$
p_n^j(x)=\sum_{k=0}^j p_{n,k}^j(\log x)^k.
$$
Proposition 3.2 implies that, for $m\geqq0$ and $0<x<1$,
$$
|Q_{N-1}^m(x)|\leqq C_2 x^{2N(m+2)}\biggl(\log\frac1x\biggr)^{m+1}.
\tag4.10
$$
In addition, it follows that, for $0\leqq h\leqq N-1$ and $m\geqq1$,
$$
Q_h^m(x)=Q_{N-1}^{m-1}(x)-\kappa_Nx^{2N(m+1)}\sum_{n=0}^h x^{2n}p_n^m(x)
\tag4.11
$$
and that, for $0\leqq h\leqq N-2$ and $m\geqq0$,
$$
Q_h^m(x)=Q_{N-1}^m(x)+\kappa_Nx^{2N(m+1)}\sum_{n=h+1}^{N-1} x^{2n}p_n^m(x).
\tag4.12
$$
The first lemma in this section is the following.

\medskip

\proclaim{Lemma 4.1}
For a given integer $L\geqq1$ we have
$$
\split
\frac{R_{N-2}(x)}{G(y)}=
&Q_{N-1}^L(y) R_{N(L+2)-2}(x)\\
&+\sum_{m=0}^L\sum_{h=0}^{N-1} a_{N(m+1)+h-1} Q_h^m(y) x^{N(m+1)+h-1}\\
&+\kappa_N\sum_{m=0}^L \sum_{h=0}^{N-1}
y^{2N(m+1)+2h}p_h^m(y) R_{N(m+1)+h-2}(x).
\endsplit
$$
\endproclaim

\medskip

The proof of Lemma 4.1 is deferred to Section 6.
For $t>0$ and integers $p,q\geqq0$ we set
$$
V(t;p,q)=\int_0^1 dy \int_0^\infty y^{2p+1}
\biggl(\log\frac1y\biggr)^q e^{-tuy}R_p(u^2)du.
$$
The other lemma is the asymptotic expansion of $V(t;p,q)$ for large $t$.

\medskip

\proclaim{Lemma 4.2}
Let $L\geqq p+1$ be a given integer. We have that, as $t\to\infty$,
$$
V(t;p,q)=\sum_{h=0}^q
\frac{\mu_{p,q}^h(\log t)^h}{t^{2p+2}}
+\sum_{n=p+1}^L\sum_{h=0}^q
 \frac{\theta_{p,q}^{n,h}(\log t)^h}{t^{2n+1}}
+O\mm{\frac{(\log t)^q}{t^{2L+3}}}
$$
for suitable constants $\mu_{p,q}^h$ and $\theta_{p,q}^{n,h}$.
\endproclaim

\medskip

The proof of Lemma 4.2 will be also given in Section 6.
We are now ready to compute (4.9). For $t>0$ let
$$
\split
&T_1(t)=\frac{2r^2}{\sqrt\pi}\int_0^1dy\int_0^\infty
\frac{Q_{N-1}^M(y)}{y^3}e^{-tuy}R_{N(M+2)-2}(u^2)du,\\
&T_2(t)=\frac{2r^2}{\sqrt\pi}\sum_{m=0}^M\sum_{h=0}^{N-1}
a_{N(m+1)+h-1}\int_0^1dy\int_0^\infty
\frac{Q_h^m(y)}{y^3}e^{-tuy}u^{2N(m+1)+2h-2}du,\\
&T_3(t)=\frac{2r^2\kappa_N}{\sqrt\pi}\sum_{m=0}^M\sum_{h=0}^{N-1}
\int_0^1dy\int_0^\infty y^{2N(m+1)+2h-3}p_h^m(y)\\
&\hskip7cm\times e^{-tuy}R_{N(m+1)+h-2}(u^2)du.
\endsplit
$$
Lemma 4.1 for $L=M$ implies that (4.9) is the sum of
these three integrals.

It follows from (4.11) for $h=N-1$ and $m=M+1$ that
$$
\split
T_1(t)&=\frac{2r^2}{\sqrt\pi}\int_0^1dy\int_0^\infty
\frac{Q_{N-1}^{M+1}(y)}{y^3}e^{-tuy}R_{N(M+2)-2}(u^2)du\\
&\hphantom{=}+\frac{2r^2\kappa_N}{\sqrt\pi}\sum_{k=0}^{M+1}
(-1)^k p_{0,k}^{M+1}V(t;N(M+2)-2,k)\\
&\hphantom{=}+\frac{2r^2\kappa_N}{\sqrt\pi}\sum_{n=1}^{N-1}\sum_{k=0}^{M+1}
(-1)^k p_{n,k}^{M+1}\int_0^1 dy \int_0^\infty y^{2n+2N(M+2)-3}
\biggl(\log\frac1y\biggr)^k\\
&\hphantom{+\frac{2r^2\kappa_N}{\sqrt\pi}\sum_{n=1}^{N-1}\sum_{k=0}^{M+1}
(-1)^k p_{n,k}^{M+1}\int_0^1 dy\int_0^\infty}
\hskip0.7cm 
\times e^{-tuy}R_{N(M+2)-2}(u^2)dy,\\
\endsplit
$$
of which the $j$-th term is denoted by $T_1^j(t)$ for $j=1,2,3$.
In virtue of (4.3) and (4.10), we have that $|T_1^1(t)|$ is bounded by
$$
C_3 \int_0^1dy\int_0^\infty y^{2N(M+3)-3}\biggl(\log\frac1y\biggr)^{M+2}
e^{-tuy}u^{2N(M+2)-2}du,
$$
which is not larger than
$$
\frac{C_4}{t^{2N(M+2)-1}}\int_0^1 y^{2N-2}\biggl(\log\frac1y\biggr)^{M+2} dy.
$$
Here we have applied (4.5). Therefore $T_1^1(t)$ is of order
$1/t^{2N(M+2)-1}$. Lemma 4.2 gives that, for $t>1$
$$
\sup_{0\leqq k\leqq M+1} |V(t;N(M+2)-2,k)|
\leqq \frac{C_5(\log t)^{M+1}}{t^{2N(M+2)-2}}.
$$
This immediately implies that  $T_1^2(t)$ is of order
$(\log t)^{M+1}/t^{2N(M+2)-2}$.
The way of estimate of $T_1^3(t)$ is similar to $T_1^1(t)$.
Indeed, applying (4.3) and (4.5), we obtain
$$
\split
|T_1^3(t)|&\leqq C_6 \sum_{n=1}^{N-1}\sum_{k=0}^{M+1}
\int_0^1 dy \int_0^\infty y^{2n+2N(M+2)-3}\biggl(\log\frac1y\biggr)^k
e^{-tuy}u^{2N(M+2)-2}du\\
&\leqq C_7 \int_0^1 dy \int_0^\infty y^{2N(M+2)-1}
\biggl(\log\frac1y\biggr)^{M+1}
e^{-tuy}u^{2N(M+2)-2}du,\\
\endsplit
$$
which is equal to a constant multiple of $1/t^{2N(M+2)-1}$.
Therefore we can conclude
$$
T_1(t)=O\mm{\frac{(\log t)^{M+1}}{t^{2N(M+2)-2}}}.
\tag4.13
$$

Carrying out the integral on $u$ in $T_2(t)$, we have by (4.5) that
$$
T_2(t)=\frac{2r^2}{\sqrt\pi}\sum_{m=0}^M\sum_{h=0}^{N-1}
\frac{b_{N(m+1)+h-1}}{t^{2N(m+1)+2h-1}}
\int_0^1\frac{Q_h^m(y)}{y^{2N(m+1)+2h+2}}dy.
\tag4.14
$$
For $0\leqq h\leqq N-2$ and $m\geqq0$ we have by (4.12) that
the integral on $y$ in the right hand side of (4.14) is the sum of
$$
\align
&\int_0^1\frac{Q_{N-1}^m(y)}{y^{2N(m+1)+2h+2}}dy,\tag4.15\\
&\kappa_N\sum_{n=h+1}^{N-1}\sum_{k=0}^m(-1)^kp_{n,k}^m
\int_0^1 y^{2n-2h-2}\biggl(\log\frac1y\biggr)^kdy.\tag4.16
\endalign
$$
It follows from (4.5) and (4.10) that the absolute value of (4.15)
is bounded by a constant multiple of 
$$
\int_0^1y^{2N-2h-2}\biggl(\log\frac1y\biggr)^{m+1}dy
=\frac{\varGamma(m+2)}{(2N-2h-1)^{m+2}}.
$$
It is easy to see that a bound of the absolute value of (4.16) is
a constant multiple of
$$
\int_0^1y^{2n-2h-2}\biggl(\log\frac1y\biggr)^m dy
=\frac{\varGamma(m+1)}{(2n-2h-1)^{m+1}}.
$$
Therefore we conclude that $Q_h^m(y)/y^{2N(m+1)+2h+2}$ is integrable
on $(0,1)$. Setting
$$
\varXi_{m,h}^2=\varXi_{m,h}^1\int_0^1\frac{Q_h^m(y)}{y^{2N(m+1)+2h+2}}dy
$$
for $0\leqq m\leqq M$ and $0\leqq h\leqq N-1$, we have
$$
T_2(t)=\frac{2r^2}{\sqrt\pi}\sum_{m=0}^M\sum_{h=0}^{N-1}
\frac{\varXi_{m,h}^1}{t^{2N(m+1)+2h-1}}.
\tag4.17
$$

The definition of $p_h^m(x)$ gives
$$
T_3(t)=\frac{2r^2\kappa_N}{\sqrt\pi}\sum_{m=0}^M\sum_{h=0}^{M-1}
\sum_{k=0}^m (-1)^kp_{h,k}^m V(t;N(m+1)+h-2,k).
$$
Lemma 4.2 for $L=N(M+2)-2$ yields that the main part of $T_3(t)$ is
the sum of
$$
\align
&\frac{2r^2}{\sqrt\pi}\sum_{m=0}^M\sum_{h=0}^{M-1}
\sum_{k=0}^m \sum_{j=0}^k 
\frac{\varXi_{m,h,k,j}^2(\log t)^j}{t^{2N(m+1)+2h-2}},\tag4.18\\
&\frac{2r^2}{\sqrt\pi}\sum_{m=0}^M\sum_{h=0}^{M-1}
\sum_{k=0}^m \sum_{n=N(m+1)+h-1}^{N(M+2)-2}\sum_{j=0}^k
\frac{\varXi_{m,h,k,n,j}^3(\log t)^j}{t^{2n+1}},\tag4.19\\
\endalign
$$
and the remainder of $T_3(t)$ is of order $(\log t)^M/t^{2N(M+2)-1}$, where
$$
\split
&\varXi_{m,h,k,j}^2=(-1)^k \kappa_N p_{h,k}^m \mu_{N(m+1)+h-2,k}^j,\\
&\varXi_{m,h,k,n,j}^3=(-1)^k \kappa_N p_{h,k}^m
\theta_{N(m+1)+h-2,k}^{n,j}.
\endsplit
$$
In addition, we put
$$
\varXi_{m,h,j}^4=\sum_{k=j}^m \varXi_{m,h,j,k}^2,\quad
\varXi_{m,h,n,j}^5=\sum_{k=j}^m \varXi_{m,h,k,n,j}^3.
$$
Then we have that (4.18) is
$$
\frac{2r^2}{\sqrt\pi}\sum_{m=0}^M\sum_{h=0}^{N-1}\sum_{j=0}^m
\frac{\varXi_{m,h,j}^4(\log t)^j}{t^{2N(m+1)+2h-2}}
\tag4.20
$$
and that (4.19) is equal to
$$
\frac{2r^2}{\sqrt\pi}\sum_{m=0}^M\sum_{h=0}^{N-1}
\sum_{n=N(m+1)+h-1}^{N(M+2)-2}\sum_{j=0}^m
\frac{\varXi_{m,h,n,j}^5 (\log t)^j}{t^{2n+1}},
$$
which can be represented by
$$
\frac{2r^2}{\sqrt\pi}\sum_{m=0}^M\sum_{h=0}^{N-1} \sum_{n=h}^{N(M-m+1)-1}
\sum_{j=0}^m \frac{\varXi_{m,h,n+N(m+1)-1,j}^5 (\log t)^j}{t^{2N(m+1)+2n-1}}.
\tag4.21
$$
We divide the summation on $n$ in (4.21)
into two parts. One is the sum over $[h,N-1]$ and the other is
the sum over $[N,N(M-m+1)-1]$. Changing the order of the summations
on $n$ and $h$ in the both sums, we obtain that (4.21) is
$$
\split
&\frac{2r^2}{\sqrt\pi}\sum_{m=0}^M\sum_{n=0}^{N-1} \sum_{j=0}^m
\biggl(\sum_{h=0}^n \varXi_{m,h,n+N(m+1)-1,j}^5\biggr)
\frac{(\log t)^j}{t^{2N(m+1)+2n-3}}\\
&\hphantom{-}+\frac{2r^2}{\sqrt\pi}\sum_{m=0}^M\sum_{n=N}^{N(M-m+1)-1}
\sum_{j=0}^m\biggl(\sum_{h=0}^{N-1} \varXi_{m,h,n+N(m+1)-1,j}^5\biggr)
\frac{(\log t)^j}{t^{2N(m+1)+2n-3}}.\\
\endsplit
$$
Let
$$
\varXi_{m,n,j}^6=\sum_{h=0}^{\min\{n,N-1\}}\varXi_{m,h,n+N(m+1)-1,j}^5
$$
and then (4.21) is represented by
$$
\frac{2r^2}{\sqrt\pi}\sum_{m=0}^M\sum_{n=0}^{N(M-m+1)-1} \sum_{j=0}^m
\frac{\varXi_{m,n,j}^6 (\log t)^j}{t^{2N(m+1)+2n-1}}.
\tag4.22
$$
Hence, with the help of (4.13), we can accordingly obtain that
(4.9) is the sum of (4.17), (4.20), (4.22) and the term
which is of order $(\log t)^{M+1}/t^{2N(M+2)-2}$.
Since the contribution of the sum on $n$ over
$[N(M-m)+1,N(M-m+1)-1]$ in (4.22)
is of order $(\log t)^M/t^{2N(M+1)+1}$, we have
$$
\split
T(t)=
&\frac{2r^2}{\sqrt\pi}\sum_{n=0}^{N-2}
\frac{\varXi_n^7}{t^{2n+1}}
+\frac{2r^2}{\sqrt\pi}\sum_{n=N-1}^{N(M+1)-2}
\frac{\varXi_n^7}{t^{2n+1}}\\
&+\frac{2r^2}{\sqrt\pi}\sum_{m=0}^M\sum_{h=0}^{N-1}
\frac{\varXi_{m,h}^1}{t^{2N(m+1)+2h-1}}\\
&+\frac{2r^2}{\sqrt\pi}\sum_{m=0}^M\sum_{h=0}^{N-1}\sum_{j=0}^m
\frac{\varXi_{m,h,j}^4(\log t)^j}{t^{2N(m+1)+2h-2}}\\
&+\frac{2r^2}{\sqrt\pi}\sum_{m=0}^M\sum_{n=0}^{N(M-m)} \sum_{j=0}^m
\frac{\varXi_{m,n,j}^6 (\log t)^j}{t^{2N(m+1)+2n-1}}\\
&+O\mm{\frac1{t^{2N(M+1)-1}}},
\endsplit
\tag4.23
$$
where
$$
\varXi_n^7=
\cases
\dis b_n \int_0^\infty\frac{dy}{y^{2n+4}G(y)}
\quad&\text{if $0\leqq n\leqq N-2$,}\\
\dis b_n \int_1^\infty\frac{dy}{y^{2n+4}G(y)}
\quad&\text{if $n\geqq N-1$.}
\endcases
$$
We remark that the second tem of the right hand side of (4.23)
can be expressed by
$$
\frac{2r^2}{\sqrt\pi}\sum_{m=0}^{M-1}\sum_{n=N(m+1)-1}^{N(m+2)-2}
\frac{\varXi_n^7}{t^{2n+1}}
=\frac{2r^2}{\sqrt\pi}\sum_{m=0}^{M-1}\sum_{n=0}^{N-1}
\frac{\varXi_{n+N(m+1)-1}^7}{t^{2N(m+1)+2n-1}}
\tag4.24
$$
and that the third term of (4.23) is
$$
\frac{2r^2}{\sqrt\pi}\sum_{m=0}^{M-1}\sum_{n=0}^{N-1}
\frac{\varXi_{m,n}^1}{t^{2N(m+1)+2n-1}}+O\mm{\frac1{t^{2N(M+1)-1}}}.
\tag4.25
$$
Since both (4.24) and the first part of (4.25) are the special cases
of the fifth term of (4.23), it follows that
$$
\split
T(t)=
&\frac{2r^2}{\sqrt\pi}\sum_{n=0}^{N-2}\frac{\varXi_n^7}{t^{2n+1}}
+\frac{2r^2}{\sqrt\pi}\sum_{m=0}^M\sum_{h=0}^{N-1}\sum_{j=0}^m
\frac{\varXi_{m,h,j}^4(\log t)^j}{t^{2N(m+1)+2h-2}}\\
&+\frac{2r^2}{\sqrt\pi}\sum_{m=0}^M\sum_{n=0}^{N(M-m)} \sum_{j=0}^m
\frac{\varXi_{m,n,j}^8 (\log t)^j}{t^{2N(m+1)+2n-1}}\\
&+O\mm{\frac1{t^{2N(M+1)-1}}}.
\endsplit
\tag4.26
$$
Here the notation $\varXi_{m,n,j}^8$ has been used to denote
$\varXi_{n+N(m+1)-1}^7+\varXi_{m,n}^1+\varXi_{m,n,0}^6$
if $0\leqq m\leqq M-1$, $0\leqq n\leqq N-1$, $j=0$ and
$\varXi_{m,n,j}^6$ otherwise.

It is easy to see that the second term of the right hand side of (4.26) is
$$
\frac{2r^2}{\sqrt\pi}\sum_{m=0}^{M-1}\sum_{n=0}^{N-1}\sum_{j=0}^m
\frac{\varXi_{m,n,j}^4(\log t)^j}{t^{2N(m+1)+2n-2}}
+\frac{2r^2}{\sqrt\pi}\sum_{j=0}^M
\frac{\varXi_{M,0,j}^4(\log t)^j}{t^{2N(M+1)-2}}
+O\mm{\frac{(\log t)^M}{t^{2N(M+1)}}}.
$$
and that the third term of (4.26) is
$$
\split
&\frac{2r^2}{\sqrt\pi}\sum_{m=0}^{M-1}\sum_{n=0}^{N(M-m)-1} \sum_{j=0}^m
\frac{\varXi_{m,n,j}^8 (\log t)^j}{t^{2N(m+1)+2n-1}}\\
&\hphantom{-}+\frac{2r^2}{\sqrt\pi}\sum_{m=0}^{M-1}\sum_{j=0}^m
\frac{\varXi_{m,N(M-m),j}^8 (\log t)^j}{t^{2N(M+1)-1}}
+\frac{2r^2}{\sqrt\pi}\sum_{j=0}^M
\frac{\varXi_{M,0,j}^8 (\log t)^j}{t^{2N(M+1)-1}},
\endsplit
$$
which is expressed by
$$
\frac{2r^2}{\sqrt\pi}\sum_{m=0}^{M-1}\sum_{n=0}^{N(M-m)-1} \sum_{j=0}^m
\frac{\varXi_{m,n,j}^8 (\log t)^j}{t^{2N(m+1)+2n-1}}
+\frac{2r^2}{\sqrt\pi}\sum_{j=0}^M
\frac{\varXi_j^9 (\log t)^j}{t^{2N(M+1)-1}}.
$$
Here we have set
$$
\varXi_j^9=\sum_{m=j}^M\varXi_{m,N(M-m),j}^8.
$$
Let
$$
W_m(t;\{\varXi_j\})=\sum_{j=0}^m \varXi_j
\biggl(\log\frac{\sqrt{2t}}r\biggr)^j
$$
for a sequence $\{\varXi_j\}_{j=0}^m$ of real numbers.
Recall that $T(t)=S(r^2t^2/2)$ for $t>0$ and
hence (4.26) implies that $S(t)$ is equal to
$$
\split
&\sum_{n=0}^{N-2}\frac{\varXi_n^7 r^{2n+3}}{2^{n-1/2}\sqrt\pi}
\frac1{t^{2n+1}}\\
&\hphantom{-}+\sum_{m=0}^{M-1}\sum_{n=0}^{N-1}
\frac{r^{2N(m+1)+2n}}{2^{N(m+1)+n-2}\sqrt\pi}
\frac{W_m(t;\{\varXi_{m,n,j}^4\})}{t^{N(m+1)+n-1}}\\
&\hphantom{-}+\sum_{m=0}^{M-1}\sum_{n=0}^{N(M-m)-1}
\frac{r^{2N(m+1)+2n+1}}{2^{N(m+1)+n-3/2}\sqrt\pi}
\frac{W_m(t;\{\varXi_{m,n,j}^8\})}{t^{N(m+1)+n-1/2}}\\
&\hphantom{-}+\frac{r^{2N(M+1)}}{2^{N(M+1)-2}\sqrt\pi}
\frac{W_M(t;\{\varXi_{M,0,j}^4\})}{t^{N(M+1)-1}}\\
&\hphantom{-}+\frac{r^{2N(M+1)+1}}{2^{N(M+1)-3/2}\sqrt\pi}
\frac{W_M(t;\{\varXi_j^9\})}{t^{N(M+1)-1/2}}\\
&\hphantom{-}+O\mm{\frac1{t^{N(M+1)-1/2}}}.
\endsplit
\tag4.27
$$
We should note that it is easy to represent $W_m(t;\{\varXi_j\})$
as the polynomial of $\log t$. Indeed, by the binomial theorem, we have
$$
W_m(t;\{\varXi_j\})
=\sum_{k=0}^m U_{m,k}(\{\varXi_j\})(\log t)^k,
\tag4.28
$$
where
$$
U_{m,k}(\{\varXi_j\})=\sum_{j=k}^m\frac{\varXi_j}{2^j}\binom jk
\biggl(\log\frac2{r^2}\biggr)^{j-k}.
$$

We must give the asymptotic expansion of the third term of
the right hand side of (2.5). The calculation is easy.
Note that $\RE(z_j^{(d)})<0$ for each $j=1,2,\dots,\phi_d$.
In virtue of (4.3) and (4.5), we can derive
$$
\int_0^\infty e^{-\frac{r^2x^2}{2t}+z_j^{(d)}x}dx
=-\sum_{n=0}^{N(M+1)-2}\frac{b_n r^{2n}}
{2^n(z_j^{(d)})^{2n+1}}\frac1{t^n}
+O\mm{\frac1{t^{N(M+1)-1}}}.
$$
This yields that the third term of (2.5) is equal to
$$
\sum_{n=0}^{N(M+1)-2} \frac{\varXi_n^{10}r^{2n+3}}{t^{n+1/2}}
+O\mm{\frac1{t^{N(M+1)-1/2}}},
\tag4.29
$$
where
$$
\varXi_n^{10}=\frac{b_n}{2^{n-1/2}\sqrt\pi}
\sum_{j=1}^{\phi_d} \frac1{(z_j^{(d)})^{2n+3}}.
$$
We split the summation on $n$ in (4.29) into two parts.
One is the sum over $[0,n-1]$ and the other is the remainder.
Thus the first term of (4.29) is
$$
\split
&\sum_{n=0}^{N-1}\frac{\varXi_n^{10}r^{2n+3}}{t^{n+1/2}}
+\sum_{m=0}^{M-1}\sum_{n=N(m+1)-1}^{N(m+2)-2}
\frac{\varXi_n^{10}r^{2n+3}}{t^{n+1/2}}\\
&=\sum_{n=0}^{N-1}\frac{\varXi_n^{10}r^{2n+3}}{t^{n+1/2}}
+\sum_{m=0}^{M-1}\sum_{n=0}^{N-1}
\frac{\varXi_{N(m+1)+n-1}^{10}r^{2N(m+1)+2n+1}}{t^{N(m+1)+n-1/2}}.
\endsplit
$$

Recall that $N=d/2-1$. We put
$$
\split
&\rho_n^{(0)}=\frac{(-1)^N\varXi_n^7}{2^{n-1/2}\sqrt\pi}+\varXi_n^{10},\quad
\rho_{m,n,k}^{(1)}=\frac{(-1)^N U_{m,k}(\{\varXi_{m,n,j}^4\})}
{2^{N(m+1)+n-2}\sqrt\pi},\\
&\rho_{m,n,k}^{(2)}=
\cases
\dis \frac{(-1)^N U_{m,0}(\{\varXi_{m,n,j}^9\})}
{2^{N(m+1)+n-3/2}\sqrt\pi}+\varXi_{N(m+1)+n-1}^{10}\quad
&\text{if $0\leqq n\leqq N-1$, $k=0$,}\\
\dis \frac{(-1)^N U_{m,k}(\{\varXi_{m,n,j}^8\})}
{2^{N(m+1)+n-3/2}\sqrt\pi}&\text{otherwise,}
\endcases\\
&\rho_k^{(3)}=\frac{(-1)^N U_{m,k}(\{\varXi_{M,0,j}^4\})}
{2^{N(m+1)-2}\sqrt\pi},\quad
\rho_k^{(4)}=\frac{(-1)^N U_{m,k}(\{\varXi_j^9\})}
{2^{N(m+1)-3/2}\sqrt\pi}.
\endsplit
$$
Hence we obtain
$$
\split
L(t)=&S_r^{(d)}\biggl[\frac{d-2}2t+\frac{r^2}{d-4}
+\sum_{n=0}^{N-2}\frac{\rho_n^{(0)}r^{2n+3}}{t^{n+1/2}}\\
&\hphantom{S_r^{(d)}\biggl[}+\sum_{m=0}^{M-1}\sum_{n=0}^{N-1}
\frac{r^{2N(m+1)+2n}}{t^{N(m+1)+n-1}}
\sum_{k=0}^m\rho_{m,n,k}^{(1)}(\log t)^k\\
&\hphantom{S_r^{(d)}\biggl[}+\sum_{m=0}^{M-1}
\sum_{n=0}^{N(M-m)-1}\frac{r^{2N(m+1)+2n+1}}{t^{N(m+1)+n-1/2}}
\sum_{k=0}^m\rho_{m,n,k}^{(2)}(\log t)^k\\
&\hphantom{S_r^{(d)}\biggl[}+\frac{r^{2N(M+1)}}{t^{N(M+1)-1}}
\sum_{k=0}^M\rho_k^{(3)}(\log t)^k
+\frac{r^{2N(M+1)+1}}{t^{N(M+1)-1/2}}\sum_{k=0}^M
\rho_k^{(4)}(\log t)^k\biggr]\\
&+O\mm{\frac1{t^{N(M+1)-1/2}}}.
\endsplit
$$
It follow from (2.6) that the first term of the smaller one
than the constant is the term of $1/t^{d/2-2}$ Since $d/2-1=N-1$,
this implies that each $\rho_n^{(0)}$ varnishes for $n=0,1,\dots,N-2$,
We finish the proof of Theorem 2.1.

We remark that $\rho_n^{(0)}=0$ yields an interesting
equality concerning zeros of a modified Bessel function.
Recall that, for $0\leqq n\leqq N-2$
$$
\varXi_n^7=b_n \int_0^\infty\frac{dy}{y^{2n+4}G(y)},\quad
\varXi_n^{10}=\frac{b_n}{2^{n-1/2}\sqrt\pi}
\sum_{j=1}^{\phi_d}\frac1{(z_j^{(d)})^{2n+3}}.
$$
Hence it follows from $\rho_n^{(0)}=0$ that
$$
\sum_{j=1}^{\phi_d}\frac1{(z_j^{(d)})^{2n+3}}
=(-1)^{N+1}\int_0^\infty\frac{dy}{y^{2n+4}G(y)}.
$$
Since $z_j^{(d)}$ is a zero of $K_N$, denoted by $z_{N,j}$, and
$\phi_d$ is the number of zeros of $K_N$, which has been
denoted by $\phi(N)$, we obtain the following.

\medskip

\proclaim{Remark 4.3}
Let $n\geqq2$ a given integer. We have that, for $m=0,1,\dots,n-2$,
$$
\sum_{j=1}^{\phi(n)}\frac1{z_{n,j}^{2m+3}}
=(-1)^{n+1}\int_0^\infty\frac{dx}{x^{2m+4}\{K_n(x)^2+\pi^2 I_n(x)^2\}}.
$$
\endproclaim

%555555555555555555555555555555555555555555555555555555

\head
5. The four dimensional case
\endhead

\noindent
This section deals with the case of $d=4$, which is equivalent to $N=1$.
In order to prove Theorem 2.2, we need to give a lemma like Lemma 4.1.
Throughout this section we use $C_8,\dots,C_{11}$ for 
positive constants independent of the variable.

Recall the definition of $p_n^j(x)$ given in Section 4.
For $m\geqq0$ and $x>0$ we set
$$
Q^m(x)=\frac1{G(x)}-x^2\sum_{j=0}^m x^{2j}p_0^j(x).
$$
Then it follows that, for $m\geqq1$ and $x>0$,
$$
Q^{m-1}(x)=Q^m(x)+x^{2m+2}p_0^m(x).
\tag5.1
$$
In addition, Proposition 3.2 implies that, for $m\geqq0$ and $0<x<1$
$$
|Q^m(x)|\leqq C_8 x^{2m+4}\biggl(\log\frac1x\biggr)^{m+1}.
\tag5.2
$$

\proclaim{Lemma 5.1}
For a given integer $L\geqq1$ we have
$$
\split
Q^0(y)e^{-x}=Q^L(y)R_{L-1}(x)
&+\sum_{m=0}^{L-1} a_m Q^m(y) x^m\\
&+\sum_{m=1}^L y^{2m+2}p_0^m(y) R_{m-1}(x).
\endsplit
$$
\endproclaim
\demo{Proof}
It follows from (5.1) that, for $m\geqq1$
$$
Q^{m-1}(y)R_{m-1}(x)=Q^m(y)R_m(x)+a_m Q^m(y) x^m
+y^{2m+2}p_0^m(y)R_{m-1}(x).
$$
Take the summation on $m$ over $[1,L]$ and then we have
$$
Q^0(y)R_0(x)=Q^L(y)R_L(x)+\sum_{m=1}^L a_m Q^m(y) x^m
+\sum_{m=1}^L y^{2m+2}p^m(y)R_{m-1}(x).
$$
The definition of $R_n(x)$ immediately yields the assertion of this lemma.
\qed
\enddemo

\pagebreak
%\medskip

We are now ready to show Theorem 2.2 and concentrate on
the second term of the right hand side of (2.4). For $t>0$ let
$$
S(t)=\frac{\sqrt 2 r^3}{\sqrt{\pi t}}
\int_0^\infty \!\!\! \int_0^\infty \frac{1-e^{-xy}}{y^3 G(y)}
e^{-\frac{r^2x^2}{2t}}dxdy
$$
and $T(t)=S(r^2t^2/2)$. A similar calculation used in Section 4 yields
$$
T(t)=\frac{2 r^2}{\sqrt\pi}\int_0^\infty dy
\int_0^\infty \frac{1-e^{-tuy}}{y^3G(y)}e^{-u^2}du.
\tag5.3
$$
Since our purpose is to give the asymptotic expansion of (5.3) for large $t$,
we may consider the case of $t>1$.
We divide the interval of the integration on $y$
into $(0,1)$ and $[1,\infty)$. We have remarked that
$1/y^3G(y)$ is integrable on $[1,\infty)$. This implies that
the integral on $y$ over this interval is easy to treat. However
we must treat the integral on $y$ over $(0,1)$ in (5.3) with much care
since (3.1) yields that $1/y^3G(y)$ is not integrable on $(0,1)$.
The right hand side of (5.3) is the sum of
$$
\align
&\frac{2r^2}{\sqrt\pi}\int_1^\infty dy\int_0^\infty
\frac{e^{-u^2}}{y^3G(y)}du,\tag5.4\\
&-\frac{2r^2}{\sqrt\pi}\sum_{n=0}^M a_n \int_1^\infty dy
\int_0^\infty \frac1{y^3G(y)}e^{-tuy}u^{2n}du,\tag5.5\\
&-\frac{2r^2}{\sqrt\pi}\int_1^\infty dy \int_0^\infty
\frac1{y^3G(y)}e^{-tuy}R_M(u^2)du,\tag5.6\\
&\frac{2r^2}{\sqrt\pi}\int_0^1dy\int_0^\infty
\frac{1-e^{-tuy}}y e^{-u^2}du,\tag5.7\\
&\frac{2r^2}{\sqrt\pi}\int_0^1 dy \int_0^\infty
\frac{1-e^{-tuy}}{y^3G(y)}Q^0(y)e^{-u^2}du.\tag5.8
\endalign
$$
It is easy to see that (5.4) and (5.5) are equal to
$$
r^2\int_1^\infty \frac{dy}{y^3G(y)},\quad
-\frac{2r^2}{\sqrt\pi}\sum_{n=0}^M b_n
\int_1^\infty \frac{dy}{y^{2n+4}G(y)}
\frac1{t^{2n+1}},
$$
respectively. Moreover we deduce from (3.1), (4.3) and (4.5)
that the absolute value of (5.6) is dominated by
$$
\frac{C_9}{t^{2M+3}}\int_1^\infty \frac{dy}{y^{2M+6}G(y)}
=O\mm{\frac1{t^{2M+3}}}.
$$
Changing the variable from $y$ to $v=ty$, we obtain
that (5.7) is
$$
\align
&\frac{2r^2}{\sqrt\pi}\int_0^1dv \int_0^\infty
\frac{1-e^{-uv}}v e^{-u^2}du\tag5.9\\
&\hphantom{-}+\frac{2r^2}{\sqrt\pi}\int_1^t dv\int_0^\infty
\frac{e^{-u^2}}v du\tag5.10\\
&\hphantom{-}-\frac{2r^2}{\sqrt\pi}\int_1^\infty dv \int_0^\infty
\frac{e^{-uv} e^{-u^2}}v du\tag5.11\\
&\hphantom{-}+\frac{2r^2}{\sqrt\pi}\int_t^\infty dv \int_0^\infty
\frac{e^{-uv} e^{-u^2}}v du.\tag5.12\\
\endalign
$$
Since $0\leqq1-e^{-uv}\leqq uv$ for $u,v\geqq0$,
the integral in (5.9) converges.
It is easy to see that (5.10) coincides with
$r^2 \log t$ and that the absolute value of (5.11) is
dominated by a constant multiple of
$$
\int_1^\infty dv \int_0^\infty \frac{e^{-uv}}v du
=\int_1^\infty \frac{dv}{v^2}=1,
$$
which implies that the integral in $(5.11)$ converges.
It follows that (5.12) is equal to
$$
\frac{2r^2}{\sqrt\pi}\sum_{n=0}^M a_n
\int_t^\infty dv \int_0^\infty \frac{u^{2n}e^{-uv}}v du
+\frac{2r^2}{\sqrt\pi}\int_t^\infty dv \int_0^\infty
\frac{e^{-uv}}v R_M(u^2)du.
\tag5.13
$$
By (4.3) and (4.5) we have that the first term of (5.13) is
$$
\frac{2r^2}{\sqrt\pi}\sum_{n=0}^M \frac{b_n}{2n+1}\frac1{t^{2n+1}}
\tag5.14
$$
and that the absolute value of the second term of (5.13) is bounded by
$$
C_{10}\int_t^\infty dv \int_0^\infty \frac{u^{2M+2}e^{-uv}}v du
=O\mm{\frac1{t^{2M+3}}}.
$$
Therefore we accordingly obtain that (5.7) is the sum of
$r^2\log t$, (5.9), (5.11), (5.14) and the part which is of order
$1/t^{2M+3}$.

We next consider (5.8). Lemma 5.1 yields that (5.8) is equal to
$$
\align
&\frac{2r^2}{\sqrt\pi}\int_0^1dy \int_0^\infty \frac{Q^0(y)}ye^{-u^2}du
\tag5.15\\
&\hphantom{-}-
\frac{2r^2}{\sqrt\pi}\int_0^1dy \int_0^\infty
\frac{Q^{M+2}(y)}{y^3}e^{-tuy}R_{M+1}(u^2) du
\tag5.16\\
&\hphantom{-}-
\frac{2r^2}{\sqrt\pi}\sum_{m=0}^{M+1} a_m
\int_0^1dy \int_0^\infty \frac{Q^m(y)}{y^3}e^{-tuy}u^{2m}du
\tag5.17\\
&\hphantom{-}-
\frac{2r^2}{\sqrt\pi}\sum_{m=1}^{M+2} \int_0^1dy \int_0^\infty
y^{2m-1}p_0^m(y)e^{-tuy}R_{m-1}(u^2)du.
\tag5.18
\endalign
$$
Calculations of all terms are easy. It follows that (5.15) is
$$
r^2\int_0^1\frac{Q^0(y)}{y^3}dy.
\tag5.19
$$
Since $|Q^0(y)|/y^3\leqq C_8 y\log(1/y)$ for $0<y<1$ obtained by (5.2),
we have that (5.19) and also (5.15)
converge. It follows from (4.3) and (5.2) that the absolute value of
(5.16) is bounded by
$$
C_{11}\int_0^1 dy \int_0^\infty y^{2M+5}\biggl(\log\frac1y\biggr)^{M+3}
u^{2M+4}e^{-tuy}du,
$$
which is not larger than a constant multiple of
$$
\frac1{t^{2M+5}}\int_0^1\biggl(\log\frac1y\biggr)^{M+3}dy
=\frac{\varGamma(M+4)}{t^{2M+5}}.
$$
Applying (4.5) again, we deduce that (5.17) is
$$
-\frac{2r^2}{\sqrt\pi}\sum_{m=0}^{M+1} b_n
\int_0^1 \frac{Q^m(y)}{y^{2m+4}} dy \frac1{t^{2m+1}}.
\tag5.20
$$
We must check that $Q^m(y)/y^{2m+4}$ is integrable on $(0,1)$.
However it can be easily shown by (5.2). Recall the definition of
$V(t;p,q)$ given in Section 4 and then we have that (5.18) is
$$
-\frac{2r^2}{\sqrt\pi}\sum_{m=1}^{M+2}\sum_{k=0}^m
(-1)^kp_{0,k}^mV(t;m-1,k).
\tag5.21
$$
Lemma 4.2 yields that the main part of (5.21) is
$$
\frac{2r^2}{\sqrt\pi}\sum_{m=1}^{M+2} \sum_{h=0}^k
\frac{\varXi_{m,h}^{11}(\log t)^h}{t^{2m}}
+\frac{2r^2}{\sqrt\pi}\sum_{n=1}^{M+2} \sum_{m=1}^n \sum_{h=0}^m
\frac{\varXi_{n,m,h}^{12}(\log t)^h}{t^{2n+1}}
\tag5.22
$$
and the remainder of (5.21) is of order $(\log t)^{M+2}/t^{2M+7}$,
where
$$
\varXi_{m,h}^{11}=\sum_{k=h}^m (-1)^{k+1} p_{0,k}^m \mu_{m-1,k}^h,\quad
\varXi_{n,m,h}^{12}=\sum_{k=h}^m (-1)^{k+1} p_{0,k}^m \theta_{m-1,k}^{n,h}.
$$
Moreover the second term of (5.22) is
$$
\split
&\frac{2r^2}{\sqrt\pi}\sum_{n=1}^{M+2} \sum_{h=1}^m \sum_{m=h}^n
\frac{\varXi_{n,m,h}^{12}(\log t)^h}{t^{2n+1}}
+\frac{2r^2}{\sqrt\pi}\sum_{n=1}^{M+2} \sum_{m=1}^n
\frac{\varXi_{n,m,0}^{12}}{t^{2n+1}}\\
&=\frac{2r^2}{\sqrt\pi}\sum_{n=1}^{M+2} \sum_{h=0}^m
\frac{\varXi_{n,h}^{13}(\log t)^h}{t^{2n+1}},
\endsplit
\tag5.23
$$
where we have set
$$
\varXi_{n,h}^{13}=\sum_{m=\max\{h,1\}}^n\varXi_{n,m,h}^{12}.
$$
This means that the leading part of (5.18) is the sum of
the first term of (5.22) and the right hand side of (5.23)
and that the remaining part of (5.18) is
of order $1/t^{2M+7}$. Therefore we can conclude that, as $t\to\infty$,
$$
\split
T(t)=
&r^2\log t+r^2\varXi^{14}
+\frac{2r^2}{\sqrt\pi}\sum_{n=1}^{M+1} \sum_{j=0}^n
\frac{\varXi_{n,j}^{11}(\log t)^j}{t^{2n}}\\
&+\frac{2r^2}{\sqrt\pi}\sum_{n=0}^M \frac{\varXi_n^{15}}{t^{2n+1}}
+\frac{2r^2}{\sqrt\pi}\sum_{n=1}^M \sum_{j=0}^n
\frac{\varXi_{n,j}^{13}(\log t)^j}{t^{2n+1}}
+O\mm{\frac{(\log t)^{M+1}}{t^{2M+2}}},
\endsplit
$$
where
$$
\split
&\varXi^{14}=\int_1^\infty\frac{dy}{y^3G(y)}
+\frac2{\sqrt\pi}\int_0^1 dv \int_0^\infty \frac{1-e^{-uv}}v e^{-u^2}du\\
&\hphantom{\varXi^{14}=}-\frac2{\sqrt\pi}\int_1^\infty dv \int_0^\infty
\frac{e^{-uv} e^{-u^2}}v du+\int_0^1\frac{Q^0(y)}{y^3}dy,\\
&\varXi_n^{15}=b_n \biggl[\frac1{2n+1}-\int_1^\infty\frac{dy}{y^{2n+4}G(y)}
-\int_0^1\frac{Q^n(y)}{y^{2n+4}}dy\biggr].
\endsplit
$$
In addition, for $n\geqq1$ let
$$
\varXi_{n,j}^{16}=
\cases
\varXi_{n,0}^{11}+\varXi_n^{15}\quad&\text{if $j=0$,}\\
\varXi_{n,j}^{11}&\text{if $j\geqq1$.}
\endcases
$$
It follows that
$$
\split
T(t)=
r^2\log t&+r^2\varXi^{14}+\frac{2r^2\varXi_0^{15}}{\sqrt\pi\, t}
+\frac{2r^2}{\sqrt\pi}\sum_{n=1}^M \sum_{j=0}^n
\frac{\varXi_{n,j}^{16}(\log t)^j}{t^{2n}}\\
&+\frac{2r^2}{\sqrt\pi}\sum_{n=1}^M \sum_{j=0}^n
\frac{\varXi_{n,j}^{13}(\log t)^j}{t^{2n+1}}
+O\mm{\frac{(\log t)^{M+1}}{t^{2M+2}}}.
\endsplit
$$

Recall that $T(t)$ is defined by $S(r^2t^2/2)$ for $t>0$.
Then we obtain
$$
\split
S(t)=
&r^2\log\frac{\sqrt{2t}}r+r^2\varXi^{14}
+\frac{2r^3\varXi_0^{15}}{\sqrt{2\pi t}}
+\sum_{n=1}^M \frac{r^{2n+2}}{2^{n-1}\sqrt\pi}
\frac{W_n(t;\{\varXi_{n,j}^{16}\})}{t^n}\\
&+\sum_{n=1}^M \frac{r^{2n+3}}{2^{n-1/2}\sqrt\pi}
\frac{W_n(t;\{\varXi_{n,j}^{13}\})}{t^{n+1/2}}
+O\mm{\frac{(\log t)^{M+1}}{t^{M+1}}}.
\endsplit
$$
The remainder of the calculation is the same as (4.27).
It follows from (2.6) that the largest order which is smaller than
the constant is $\log t/t$. This immediately implies that
$\varXi_0^{15}=0$. We put
$$
\sigma^{(1)}=\varXi^{14}+\frac12 \log\frac2{r^2},\quad
\sigma_{n,k}^{(2)}=\frac{U_{n,k}(\{\varXi_{n,j}^{16}\})}{2^{n-1}\sqrt\pi},\quad
\sigma_{n,k}^{(3)}=\frac{U_{n,k}(\{\varXi_{n,j}^{13}\})}{2^{n-1/2}\sqrt\pi}
$$
and hence obtain by (4.28) that
$$
\split
S(t)=&
\frac12r^2\log t+r^2\sigma^{(1)}+
\sum_{n=1}^M\frac{r^{2n+2}}{t^n}\sum_{k=0}^n
\sigma_{n,k}^{(2)}(\log t)^k\\
&+\sum_{n=1}^M\frac{r^{2n+3}}{t^{n+1/2}}\sum_{k=0}^n
\sigma_{n,k}^{(2)}(\log t)^k
+O\mm{\frac{(\log t)^{M+1}}{t^{M+1}}}.
\endsplit
$$
This completes the proof of Theorem 2.2.

Similarly to Remark 4.3, we can easily deduce the following equality
from the fact that $\varXi_0^{17}$ varnishes.

\medskip

\proclaim{Remark 5.1}
We have that
$$
\int_0^\infty \frac1{x^4}
\biggl[\frac1{K_1(x)^2+\pi^2 I_1(x)^2}-x^2\biggr]dx=0.
$$
\endproclaim

%6666666666666666666666666666

\head
6. The proofs of Lemmas.
\endhead

\noindent
This section is devoted to proving lemmas given in Section 4.
We first try to see Lemma 4.1. Let $N\geqq2$ and $M\geqq1$ be integers.
We note that, for $1\leqq h\leqq N-1$ and $m\geqq0$,
$$
Q_{h-1}^m(y)=Q_h^m(y)+\kappa_Ny^{2N(m+1)+2h}p_h^m(y)
$$
and then obtain
$$
\split
Q_{h-1}^m(y)R_{N(m+1)+h-2}(x)=
&Q_h^m(y)R_{N(m+1)+h-1}(x)\\
&+a_{N(m+1)+h-1}Q_h^m(y) x^{N(m+1)+h-1}\\
&+\kappa_Ny^{2N(m+1)+2h}p_h^m(y)R_{N(m+1)+h-2}(x).
\endsplit
$$
This immediately yields that, for $m\geqq0$,
$$
\split
&Q_0^m(y)R_{N(m+1)-1}(x)\\
&=Q_{N-1}^m(y)R_{N(m+2)-2}(x)\\
&\hphantom{=}+\sum_{h=1}^{N-1} a_{N(m+1)+h-1}Q_h^m(y)x^{N(m+1)+h-1}\\
&\hphantom{=}+\kappa_N \sum_{h=1}^{N-1}
y^{2N(m+1)+2h}p_h^m(y)R_{N(m+1)+h-2}(x).
\endsplit
\tag6.1
$$
It follows from (4.11) that, for $m\geqq1$,
the left hand side of (6.1) is
$$
\split
Q_{N-1}^{m-1}(y)R_{N(m+1)-2}(x)&-\kappa_N y^{2N(m+1)}p_0^m(y)R_{N(m+1)-2}(x)\\
&-a_{N(m+1)-1}Q_0^m(y)x^{N(m+1)-1}.
\endsplit
$$
This implies that, for $m\geqq0$,
$$
\split
Q_{N-1}^{m-1}(y)R_{N(m+1)-2}(x)=
&Q_{N-1}^m(y)R_{N(m+2)-2}(x)\\
&+\sum_{h=0}^{N-1} a_{N(m+1)+h-1}Q_h^m(y)x^{N(m+1)+h-1}\\
&+\kappa_N\sum_{h=0}^{N-1}y^{2N(m+1)+2h}p_h^m(y)R_{N(m+1)+h-2}(x).
\endsplit
$$
Taking the summation on $m$ over $[1,L]$, we can derive
$$
\split
&Q_{N-1}^0(y)R_{2N-2}(x)\\
&=Q_{N-1}^L(y)R_{N(L+2)-2}(x)\\
&\hphantom{=}+\sum_{m=1}^L\sum_{h=0}^{N-1}
a_{N(m+1)+h-1}Q_h^m(y)x^{N(m+1)+h-1}\\
&\hphantom{=}+\kappa_N \sum_{m=1}^L \sum_{h=0}^{N-1}
y^{2N(m+1)+2h}p_h^m(y)R_{N(m+1)+h-2}(x).
\endsplit
\tag6.2
$$
In the case of $m=0$, we have that (6.1) means
$$
\split
Q_0^0(y)R_{N-1}(x)=
&Q_{N-1}^0(y)R_{2N-2}(x)
+\sum_{h=1}^{N-1} a_{N+h-1}Q_h^0(y)x^{N+h-1}\\
&+\kappa_N \sum_{h=1}^{N-1} y^{2N+2h}p_h^0(y)R_{N+h-2}(x).
\endsplit
$$
Since $Q_0^0(y)R_{N-1}(x)$ is represented by
$$
\frac{R_{N-2}(x)}{G(y)}-a_{N-1} Q_0^0(y)x^{N-1}
-\kappa_N y^{2N}p_0^0(y)R_{N-2}(x),
$$
it follows that
$$
\split
\frac{R_{N-2}(x)}{G(y)}=
&Q_{N-1}^0(y)R_{2N-2}(x)
+\sum_{h=0}^{N-1} a_{N+h-1}Q_h^0(y)x^{N+h-1}\\
&+\kappa_N \sum_{h=0}^{N-1} y^{2N+2h}p_h^0(y)R_{N+h-2}(x).
\endsplit
\tag6.3
$$
Combining (6.2) with (6.3), we can conclude the assertion of Lemma 4.1.

We next show Lemma 4.2. A change of the variable from $y$ to $v=ty$
yields
$$
\split
V(t;p,q)
&=\frac1{t^{2p+2}}\int_0^t dv \int_0^\infty v^{2p+1}
\biggl(\log\frac tv\biggr)^q e^{-uv}R_p(u^2)du\\
&=\frac1{t^{2p+2}}\sum_{j=0}^q\binom qj (\log t)^j \varLambda(t;p,q-j),
\endsplit
\tag6.4
$$
where
$$
\varLambda(t;\alpha,\beta)=\int_0^t dv \int_0^\infty v^{2\alpha+1}
\biggl(\log\frac1v\biggr)^\beta e^{-uv}R_\alpha(u^2)du
$$
for integers $\alpha,\beta\geqq0$. We have that $\varLambda(t;\alpha,\beta)$
is equal to
$$
\align
&\int_0^\infty dv\int_0^\infty v^{2\alpha+1}\biggl(\log\frac1v\biggr)^\beta
e^{-uv}R_\alpha(u^2)du\tag6.5\\
&\hphantom{-}-\int_t^\infty dv\int_0^\infty v^{2\alpha+1}
\biggl(\log\frac1v\biggr)^\beta e^{-uv}R_\alpha(u^2)du.\tag6.6
\endalign
$$
In order to see that (6.5) converges, we divide the domain of
the integration into the following three parts:
$$
\text{(i) $v\geqq1$, $u>0$,}\quad
\text{(ii) $0<v<1$, $0<u<1$,}\quad
\text{(iii) $0<v<1$, $u\geqq1$.}
$$
It follows from (4.3) and (4.5) that
$$
\biggl| \int_1^\infty dv\int_0^\infty v^{2\alpha+1}
\biggl(\log\frac1v\biggr)^\beta e^{-uv}R_\alpha(u^2)du\biggr|
\leqq A_\alpha\varGamma(2\alpha+3) \int_1^\infty \frac{(\log v)^\beta}{v^2}dv.
$$
This implies that the double integral in Case (i) converges.
Using (4.3) again, we have that the absolute value
of the double integral in Case (ii) is bounded by
$$
A_\alpha\int_0^1 dv \int_0^1  v^{2\alpha+1}
\biggl(\log\frac1v\biggr)^\beta e^{-uv} u^{2\alpha+2}du.
\tag6.7
$$
We dominate $v^{2\alpha+1}$, $e^{-uv}$, $u^{2\alpha+2}$ by $1$ and apply (4.5)
to the integral on $v$. Hence a bound of (6.7) is $A_\alpha\varGamma(\beta+1)$.
Applying (4.4) and (4.5), we have
$$
\biggl| \int_0^1 dv\int_1^\infty v^{2\alpha+1}
\biggl(\log\frac1v\biggr)^\beta e^{-uv}R_\alpha(u^2)du\biggr|
\leqq B_\alpha\varGamma(2\alpha+1) \int_0^1\biggl(\log\frac1v\biggr)^\beta dv.
$$
This means that the double integral in Case (iii) converges.
Therefore we conclude that (6.5) is finite and then we write
$\bar\mu_{\alpha,\beta}$ for (6.5). Let $L\geqq\alpha+1$ be an integer.
The definition of $R_n(x)$ gives that (6.6) is the sum of
$$
\align
&-\sum_{n=\alpha+1}^L a_n \int_t^\infty dv
\int_0^\infty v^{2\alpha+1}\biggl(\log\frac1v\biggr)^\beta e^{-uv}u^{2n}du,
\tag6.8\\
&-\int_t^\infty dv \int_0^\infty
v^{2\alpha+1}\biggl(\log\frac1v\biggr)^\beta e^{-uv} R_L(u^2)du.
\tag6.9
\endalign
$$
It is easy to see that (6.9) is the error term. Indeed, it follows from
(4.3) and (4.5) that, for $t>1$, the absolute value of (6.9)
is dominated by
$$
A_L\varGamma(2L+3)\int_t^\infty \frac{(\log v)^\beta}{v^{2L-2\alpha+2}},
$$
which is of order $(\log t)^\beta/t^{2L-2\alpha+1}$.
We deduce from (4.5) that, for $t>1$, (6.8) is equal to
$$
(-1)^{\beta+1}\sum_{n=\alpha+1}^L b_n
\int_t^\infty \frac{(\log v)^\beta}{v^{2n-2\alpha}}dv.
\tag6.10
$$
Note that, for $t>1$ and $n\geqq \alpha+1$, the integral on $v$ in (6.10) is
$$
\int_{\log t}^\infty x^\beta e^{-(2n-2\alpha-1)x} dx
=\frac1{t^{2n-2\alpha-1}}\sum_{h=0}^\beta
\frac{\varGamma(\beta+1)(\log t)^h}{\varGamma(h+1)(2n-2\alpha-1)^{\beta-h+1}}
$$
(see \cite{4, p.340}). Hence we obtain that (6.10) is expressed by
$$
\sum_{n=\alpha+1}^L\sum_{h=0}^\beta
\frac{\bar\theta_{\alpha,\beta}^{n,h}(\log t)^h}{t^{2n-2\alpha-1}},
$$
where
$$
\bar\theta_{\alpha,\beta}^{n,h}
=\frac{(-1)^{\beta+1} b_n \varGamma(\beta+1)}
{\varGamma(h+1)(2n-2\alpha-1)^{\beta-h+1}}.
$$
Therefore we deduce from (6.4) that
$$
\split
V(t;p,q)=
&\sum_{j=0}^q \binom qj \frac{\bar\mu_{p,q} (\log t)^j}{t^{2p+2}}
+\sum_{n=p+1}^L \sum_{j=0}^q \sum_{h=0}^{q-j}
\binom qj \frac{\bar\theta_{p,q-j}^{n,h} (\log t)^h}{t^{2n+1}}\\
&+O\mm{\frac{(\log t)^q}{t^{2L+3}}}.
\endsplit
$$
We set
$$
\mu_{p,q}^j=\binom qj \bar\mu_{p,q},\quad
\theta_{p,q}^{n,h}=\sum_{j=0}^{q-h} \binom qj \bar\theta_{p,q-j}^{n,h}
$$
and then finish the proof of Lemma 4.2.

%RRRRRRRRRRRRRRRRRRRRRRRRRRRRRRRRRRRRRRRR

\enddocument